\documentclass[12pt]{article}
\usepackage{amsmath, amsfonts, amssymb, amsthm} 
\newtheorem{theorem}{Theorem}[section]
\newtheorem{proposition}[theorem]{Proposition}
\newtheorem{definition}[theorem]{Definition}
\newtheorem{remark}[theorem]{Remark}

\newtheorem{examples}[theorem]{Examples}

\newtheorem{corollary}[theorem]{Corollary}
\newtheorem{lemma}[theorem]{Lemma}

\def\R{\mathbb R} \def\Z{\mathbb Z} \def\C{\mathbb C} 

\def\N{\mathbb N}
\def\C{{\mathbb C}} 

\def\P{\mathbb P}

\def\<{\,<\!}
\def\>{\!>\,}
\usepackage{graphics}
\usepackage{epsf}
\usepackage{epsfig}
\begin{document}

\title{Continued fractions for complex numbers \\
and values of binary quadratic forms}

\author{S.G. Dani and Arnaldo Nogueira{}\\}

\maketitle
\begin{abstract}
\noindent
\end{abstract}
We describe various properties of continued fraction expansions of complex 
numbers in terms of Gaussian integers. Numerous distinct such expansions 
are possible for a complex number. They can be arrived at through various
algorithms, as also in a more general way from what we call
``iteration sequences''. We consider in this broader context the analogues of 
the Lagrange theorem characterizing  quadratic surds, the growth properties 
of the denominators of the convergents, and the overall relation between sequences
satisfying certain conditions, in terms of nonoccurrence of certain finite
blocks, and the sequences involved in continued fraction
expansions. The results are also applied to describe a class of binary
quadratic forms with complex coefficients whose values over the set 
of pairs of Gaussian integers form a dense set of complex numbers.  
\section{Introduction} 

Continued fraction expansions of real numbers have been a classical 
topic in Number theory, playing an important role in Diophantine
approximation and various other areas. Considerable work has also been 
done on multi-dimensional versions of the
concept. There has however been relatively little work on 
analogous continued fraction expansion of complex numbers in terms of 
Gaussian integers.  A beginning was indeed made by
Hurwitz in \cite{H} as far back as 1887, where certain properties of the
nearest integer algorithm were studied and an analogue of the
classical Lagrange theorem for quadratic surds was proved. 
Some further properties of the Hurwitz continued fraction expansions 
have been described recently by Hensley \cite{Hen}.
In \cite{L}
Leveque explored another direction, focussing primarily on the aspect
of relating continued fractions to best approximations in analogy with
the classical case. In another work 
A. Schmidt \cite{Sch}  introduces an entirely
different approach to the phenomenon of continued fraction
expansions; while the results have some relation with certain classical
results the overall framework is quite different from the classical case.

This article is an attempt to explore how far the spirit of classical
continued fractions for real numbers can be extended to complex
numbers, in terms of the Gaussian integers, viz. $x+iy$ with $x$ and $y$
integers. We view a sequence $\{a_n\}_{n=0}^\infty$, of Gaussian 
integers, with $a_n\neq 0$ for $n\geq 1$, as a continued fraction
expansion of a 
complex number $z$ if the corresponding convergents can be
formed in the usual way (see below for details) from $a_0,a_1, \dots, a_n$ 
for all $n$ in analogy with the
usual convergents and they converge to $z$. Each complex number has
uncountably many distinct continued fraction expansions in this sense. 

A class of continued fraction expansions can be generated through algorithms 
for describing the successive partial quotients, of which the Hurwitz
algorithm is an  
example. Various examples of such algorithms are described in \S 2. In \S 3 we
discuss a more general approach to continued fraction expansions of complex
numbers, introducing the notion of ``iteration sequences'', the terms
of which  
correspond to the tails of continued fraction expansions. We discuss
the asymptotics  
of these sequences in a general framework and prove the convergence of 
the ``convergents'' associated with them, which is seen to give a
comprehensive 
picture of the possibilities for continued fractions of complex numbers; see 
Theorem~\ref{prop2}.

In \S 4 we consider continued fraction expansions of quadratic surds, a 
classical topic for the usual continued fractions. We prove the analogue 
of the Lagrange theorem in the general framework, and show in particular
that for a large class of algorithms (satisfying only a 
mild condition) the continued fraction expansion of a complex number
$z$ is eventually  
periodic if and only if $z$ is  a quadratic surd; see 
Corollary \ref{lagrange}.

Sections  5 and 6 are devoted to the behaviour of the convergents associated to general 
sequences of partial quotients. In particular we describe
conditions on the sequences of partial quotients, in terms of nonoccurrence of
certain blocks, which ensure that 
the denominators of the convergents have increasing absolute 
values and exponential growth. This facilitates in setting up spaces of 
sequences and continuous maps onto the space of irrational complex numbers 
such that each sequence is a continued fraction expansion of its image under
the map; see Theorem~\ref{cgce}. We exhibit also a space of sequences described
purely in terms of nonoccurrence of finitely many blocks of length two, each of 
which serves as a continued fraction expansion of a complex number, and conversely
every complex number admits an expansion in terms of a such a sequence;
see Theorem~\ref{represent}.

In \S 7 we relate the complex continued fractions to fractional
linear transformations of the complex plane and the projective space. 
Associating to each initial block of the continued fraction expansion of 
a number $z$ a projective transformation in a canonical fashion, we show that the 
corresponding sequence of transformations asymptotically takes all points of the
projective space to a point associated with $z$; see Theorem~\ref{cfl}. 

Results from \S 6 and \S 7 are applied in \S 8 to study the question of 
values of binary quadratic forms, with complex coefficients, over the set of pairs 
of Gaussian integers being dense in complex numbers. It has been known 
that for a nondegenerate quadratic form in $n\geq 3$ complex variables 
the set of values over the $n$-tuples of Gaussian integers is dense in 
complex numbers, whenever the form is not a multiple of a form with 
rational coefficients; see \S 8 for some details on the problem. The corresponding  
statement is however not true for binary quadratic forms.  From general dynamical 
considerations the set of values can be deduced to be dense in $\C$
for ``almost 
all'' binary quadratic forms (see \S 8 for details), but no specific
examples are  
known. We describe a set of complex numbers, called ``generic'', such 
that the values of
the quadratic form $(z_1-\varphi z_2)(z_1-\psi z_2)$, where $\varphi, \psi \in \C$, over the set 
with $z_1, z_2$ Gaussian integers are dense in $\C$, if $\varphi$ or $\psi$  is generic; see 
Corollary \ref{cor:qf}.
The set of generic complex numbers contains a dense ${\cal G}_\delta$ subset of $\C$, 
and using a result of D. Hensley \cite{Hen} it can also be seen to have  full Lebesgue 
measure. 

\section{Preliminaries on continued fraction expansions}

Let $\frak G$ denote the ring of Gaussian integers in $\C$, viz. $\frak
G=\{x+iy\mid x, y\in \Z\}$.  Let $\{a_n\}_{n=0}^\infty$ be a sequence
in $\frak G$. We associate to it two
sequences $\{p_n\}_{n=-1}^\infty$ and $\{q_n\}_{n=-1}^\infty$  defined 
recursively by the relations
$$ p_{-1}=1, p_0=a_0, p_{n+1}=a_{n+1}p_n+p_{n-1}, \; \mbox{for all }
n\geq 0, 
$$
$$
\mbox{ and }\ \  q_{-1}=0, q_0=1, q_{n+1}=a_{n+1}q_n+q_{n-1}, \; \mbox{
  for all } n \geq 0.
$$
We call $\{p_n\}, \{q_n\}$ the {\it $\cal Q$-pair of sequences associated
to $\{a_n\}$}; ($\cal Q$ signifies ``quotient'').  We note that 
$p_nq_{n-1}-q_np_{n-1}=(-1)^{n-1}$ for all $n\geq 1$.
If $q_n\neq 0$ for all $n$ then we can form the ``convergents''
$p_n/q_n$ and if they converge, as $n\to \infty$ to a complex number
$z$ then $\{a_n\}$ can be thought of as a continued fraction expansion
of $z$; when this holds we say that {\it $\{a_n\}$ defines a continued
  fraction expansion of $z$}, and write $z=\omega (\{a_n\})$. We note
that this confirms  
to the usual idea of continued fraction expansion as $z=a_0+
\displaystyle{\frac 1{a_1+\frac 1{a_2+\cdots}}}$ (see \cite{HW} for example).

Also, given a finite sequence $\{a_n\}_{n=0}^{m}$, where $m\in \N$
we associate to it sequences $\{p_n\}_{n=-1}^m$ and $\{q_n\}_{n=-1}^m$ 
by the recurrence relations as above, and it may also be called the
corresponding 
$\cal Q$-pair; in this case, provided all $q_n$ are nonzero we set 
$\omega (\{a_n\})=p_m/q_m$.

Though there will be some occasions to deal with continued fraction 
expansions of rational complex numbers (namely $z=x+iy$ with 
$x$ and $y$ rational), by and large we shall be concerned 
with continued fraction expansions of complex numbers which are not
rational (namely $z=x+iy$ with either $x$ or $y$ irrational);
this set of complex numbers will be denoted by $\C'$. 

We note here the following general fact about the convergence 
of the sequence of the ``convergents''.

\begin{proposition}\label{prop:cgcerate}
Let $\{a_n\}$ be a sequence  in $\frak G$ defining a continued
fraction expansion of a $z\in \C'$ and $\{p_n\}$, $\{q_n\}$ be 
the corresponding $\cal Q$-pair. Suppose that there exist $r_0\in \N$ 
and $\theta >1$ such that $|q_{n+r}|\geq \theta |q_{n}|$
for all $n\geq 0$ and $r\geq r_0$. Then there exists a $C>0$ 
such that $|z-\frac{p_{n}}{q_{n}}|\leq C |q_{n}|^{-2}$ for all $n\in \N$. 

\end{proposition}

\proof  We have $p_{n+1}q_{n}-p_{n}q_{n+1}=(-1)^{n}$ for all
$n\geq 0$, so
$\displaystyle \frac {p_{n+1}}{q_{n+1}}-\frac {p_{n}}{q_{n}} =\frac
  {(-1)^{n}}{q_{n}q_{n+1}}$, and in turn 
$\displaystyle |\frac {p_{n+m}}{q_{n+m}}-\frac {p_{n}}{q_{n}}|\leq
\sum_{k=0}^{m-1} \frac   {1}{|q_{n+k}q_{n+k+1}|}$ for all $n\geq 0$ and 
$m\geq 1$. Therefore 
$${ |z-\frac {p_{n}}{q_{n}}|\leq  \sum_{k=0}^{\infty}  
\frac   {1}{|q_{n+k}q_{n+k+1}|}}\leq r_0(1-\theta^{-1})^{-1}\frac 1{|q_n|^2},$$ 
which proves the proposition, with the choice $C=
r_0(1-\theta^{-1})^{-1}$. \hfill $\Box$

\medskip
There is considerable flexibility in writing
continued fraction expansions of complex numbers, as will become apparent
as we proceed. One of the simplest ways of arriving at such 
expansions is via algorithms for generating a continued fraction expansion
starting with a given complex number.  Through the rest of this section we shall
discuss various possibilities in this respect, together with their
significance.  

We begin by introducing some notation. 
For any $z\in \C$ and $r>0$ we shall denote by $B(z,r)$ the open disc with
center at $z$ and radius $r$, namely $\{\zeta \in \C\mid |\zeta-z|<r\}$, and
by $\bar B(z,r)$ its closure $\{\zeta \in \C\mid |\zeta -z|\leq
r\}$. Furthermore, we shall also write, for $z\in \C$, $B(z)$ for
$B(z,1)$,  $\bar B (z)$ for $\bar B(z,1)$, $B$ for $B(0,1)$  and $\bar B$
for $\bar B(0,1)$. The set of all 
nonzero complex numbers will be denoted by $\C^*$.

Let  $f:\C^* \to \frak G $ be a map such
that $f(z)\in \bar B(z)$ for all $z\in \C^*$.  There are multiple choices 
possible for the values of $f$ at any point. 
Using the map we associate to each $z\in \C$ two sequences $\{z_n\}$ and 
$\{a_n\}$ (finite or infinite) as follows. 
We define $z_0=z$, $a_0=f(z)$ and having defined 
$z_0, \dots , z_{n}$ and  $a_0, \dots , a_{n}$, $n\geq 0$, we terminate the 
sequences if $f(z_n)=z_n$ and
otherwise define $z_{n+1}=(z_n-f(z_n))^{-1}$ and $a_{n+1}=f(z_{n+1})$.
The map $f$ will  be called the {\it choice function} of the algorithm, and 
$\{a_n\}$ will be called the {\it continued fraction expansion of}
$z$ with respect to the algorithm. 
When  the sequences are terminated as above, say with $z_m$ and
$a_m$ where $m\in \N$, then it can be verified inductively that 
$z=p_m/q_m$, where $\{p_n\}$, $\{q_n\}$ is the corresponding $\cal Q$-pair. 
The crucial case will be when the sequences are infinite, which will necessarily 
be the case when $z$ is irrational. The following theorem is a particular case of 
Theorem~\ref{prop2}. 

\begin{theorem}\label{algcgce}
Let $f:\C^* \to \frak G $ be the choice function of an algorithm such  
that  $|f(z)-z|\neq 1$ for any $z$ with $|z|=1$.  Let  $z\in \C'$ be
given and $\{a_n\}$ be the  
sequence of partial quotients associated to it with respect to the algorithm 
corresponding to $f$. Let $\{p_n\}$, $\{q_n\}$ be the corresponding
$\cal Q$-pair.  
Then $q_n\neq 0$ for all $n$ and $\{p_n/q_n\}$ converges to $z$.
\end{theorem}

 This condition in the theorem is satisfied in particular if we 
assume at the outset that $f(z)\in  B(z)$, in place of $\bar B(z)$, for all $z$, but retaining 
the generality seems worthwhile.

In some contexts it is convenient to have a slightly more general
notion in which the choice function may be allowed to be multi-valued
over a (typically small, lower dimensional) subset of $\C$, and the algorithm proceeds 
by picking one
of the values. Thus we pick a map $\varphi : \C'\to {\cal P}(\frak G)$ where 
${\cal P}(\frak G)$ is the set of subsets of $\frak G$, such 
that $\varphi (z)\subset   \bar B (z)$ for all 
$z\in \C'$; we choose  $z_0=z$, $a_0\in \varphi (z)$ and having chosen 
$z_0, \dots , z_{n}$ and  $a_0, \dots , a_{n}$, $n\geq 0$,
we terminate the sequences if $z_n\in \varphi (z_n)$ and otherwise
define $z_{n+1}=(z_n-a_n)^{-1}$ and $a_{n+1}\in 
\varphi (z_{n+1})$. We call $\varphi$ as above a {\it multivalued choice
  function} for an algorithm.   One
advantage of considering multivalued choice functions $\varphi$ is
that these can be chosen to have symmetry properties like $\varphi
(\overline z)=\varphi (z)$ and $\varphi (-z)=\varphi (z)$ for all
$z$. Given a multivalued choice function one can get a
(single-valued) choice function by auxiliary conventions; in
general it may not be possible to ensure the latter to have the
symmetry properties. 

For the algorithm defined by a choice function $f$  (respectively by
a multivalued choice function $\varphi$) we call  the subset 
$\overline {\{z-f(z)\mid z\in \C'\}}$ (respectively $\overline {\{z-\alpha 
\mid \alpha \in \varphi (z), z\in  \C'\}}$) the {\it fundamental set} of the
algorithm, and denote  it by $\Phi_f$ (respectively $\Phi_\varphi$), or simply
$\Phi$. The fundamental sets play a crucial role in the 
properties of the associated algorithm. 

\begin{examples}
{\rm 
\ 

{\bf 1}. One of the well-known
algorithms is the {\it nearest integer algorithm}, discussed in detail 
by Hurwitz \cite{H}, known also as the Hurwitz algorithm. It consists of
assigning to $z\in \C$ the Gaussian integer nearest to it; this may
be viewed as a multivalued choice function $\varphi :\C\to {\cal
  P}(\frak G)$ defined by $\varphi (z) =\{\alpha \in \frak G\mid
|z-\alpha| \leq |z- \gamma| \mbox { for all } \gamma \in \frak G\}$.  
 In this case $\Phi$ is the square $\{x+iy\mid |x|\leq \frac 1{2}, |y|\leq \frac
1{2}\}$. 

{\bf 2}. The nearest integer algorithm may be generalized to a class
of algorithms as follows: 
Let $\psi :\C^* \to (0,\infty)$ be a function, and $r\in [\frac
1{\sqrt 2}, 1]$ be given. We define a multivalued choice function $\varphi$
by choosing,  
for $z\in \C'$, $\varphi (z)$ to be the set of points $x$ in $\bar
B(z,r)\cap \frak G$
where $\psi (x-z)$ is minimum. The nearest integer algorithm is a
special case with $\psi (z)=|z|$ for all $z\in \C^*$. A class of 
algorithms $f_d$, where $d\in \C$, may be defined by considering the functions 
defined by $\psi (z)=|z-d|$, in place of $|z|$; we call $\varphi_d$
 the {\it shifted Hurwitz algorithm} with displacement $d$. We note that 
 for the algorithm $\varphi_{1/2}$ the continued fraction 
 expansions of real numbers coincide with the usual simple continued 
 fraction expansions.  It may also be seen that
for $d$ with $|d|< 1-\frac 1{\sqrt 2}$, $\Phi_{\varphi_d}$ is contained in 
$\bar B(0,|d|+\frac 1{\sqrt 2})$.  

{\bf 3}. Let $f:\C^* \to \frak G$ be the map defined as follows: let $z\in \C^*$ and 
$z_0\in \frak G$ be such that $z-z_0=x+iy$ with $0\leq x <1$ and 
$0\leq y< 1$; if $x^2+y^2<1$ we define $f(z)$ to be $z_0$ and otherwise 
$f(z)$ is defined to be $z_0+1$ or $z_0+i$ whichever is nearer to $z$ (choosing
say the former if they are equidistant). It may be seen that with respect to  this algorithm
if $z\in \C$ and $\{a_n\}$ is the corresponding sequence of partial quotients 
then for all $n\geq 1$ we have $\hbox{\rm Re } a_n\geq -1$ and $\hbox{\rm Im } a_n\geq -1$;
this may be compared with the simple continued fractions for real numbers where 
the later partial quotients are positive. 

{\bf 4}.  Let $\mu :\frak G\to
(0,\infty)$  be a 
function on the set of Gaussian integers.  Let $\frac 1{\sqrt 2} \leq
r\leq 1$. A multivalued choice function $\varphi$ for an algorithm may
be defined 
by setting $\varphi (z)$, for $z\in \C'$ to be the subset of  $\bar
B(z,r)\cap \frak G$ consisting of points where the value of $\mu$ 
over the set is maximum. 

{\bf 5}. One may also restrict the range of the choice function to a
subset of  $\frak G$. If $\frak G'$ is a subset of $\frak G$ such that
for some $r\leq 1$, $B(z,r) \cap \frak G'$ is nonempty for all $z\in
\C'$ then choice functions $f$ can be constructed with values in
$\frak G'$, by adopting any of the above strategies. Here is an
example. We shall say that $x+iy \in \frak G$ is {\it even} if $x+y$ is
even. Let $\frak G'$ be the set of all even elements in $\frak G$; it
is the subgroup generated by $1+i$ and $1-i$, and one can see that 
$B(z,r) \cap \frak G'$ is nonempty for all $z\in \C'$. In this case by
analogy with the nearest integer algorithm we get the ``nearest even
integer'' algorithm; for this algorithm $\Phi$ is the square with vertices
at $\pm 1$ and $\pm i$.  

In \S 6 we discuss an algorithm with certain interesting properties, 
involving choice of nearest 
integers as well as nearest odd integers depending on the point. 
}
\end{examples}

\begin{proposition}\label{rational}
Let $f:\C \to \frak G$ be a choice function such that $\Phi_f$ is contained
in $B(0,r)$ for some $0<r<1$. Then for any rational complex number the
continued fraction expansion with respect to $f $ terminates. 
\end{proposition}

\proof Let $z=a/b$ be a rational complex number, where $a,b\in \frak G$ and 
$b\neq 0$, and $\{z_n\}$ and $\{a_n\}$
be the sequences, finite or infinite, as above starting with $z$, with respect 
to the algorithm given by $f$. Let 
$\{p_n\}$, $\{q_n\}$ be the $\cal Q$-pair corresponding to $\{a_n\}$. 
From the recurrence relations it may be seen inductively that $q_nz-p_n=
(-1)^n(z_0-f(z_0))(z_1-f(z_1)) \cdots (z_n-f(z_n))$ for all $n\geq 0$, until it is 
defined (see the proof of Proposition~\ref{prop} below). 
The condition in the hypothesis then implies that $|q_nz-p_n|< r^{n+1}$
for all $n$ as above. For $n$ such that $r^{n+1}\leq |b|^{-1}$ the 
condition can hold only if $q_nz-p_n=0$. From the equality as above 
this implies that  $f(z_n)=z_n$ for some $n$, so the sequences terminate. 
This proves the proposition. 
\hfill $\Box$

The algorithms discussed above are  time-independent
``Markovian'' in the choice strategy, where $z_{n+1}$, $n\geq 1$, is chosen
depending only on $z_{n}$ (and independent of $z_m$, $0\leq m \leq
n-1$). For certain purposes it is also of interest to consider
``time-dependent'' algorithms. An interesting example of this is seen
in the work of Leveque \cite{L}, where an algorithm is introduced
such that the corresponding continued fraction expansions have the 
best approximation property like the classical continued fractions;
this is not shared by the nearest integer algorithm (and also perhaps
by other time-independent algorithms). We shall not specifically
deal with this aspect. Our discussion in the following sections will 
however include the continued fraction expansions arising through such algorithms 
as well. 

\section{Iteration sequences}

In this section we discuss  another, more general approach to the issue
of continued fraction expansions,
introducing what we call iteration sequences. In this respect
for convenience we restrict to irrational complex numbers. 
Let $z\in \C'$ be given. We call a sequence $\{z_n\}_{n=0}^\infty$ 
an {\it iteration sequence} for $z$ if 
$$z_0=z, z_0-z_1^{-1}\in \frak G, \hbox{\rm and for all  }n\geq 1, \ |z_n|\geq 1 \mbox { and } 
z_{n}-z_{n+1}^{-1} \in \frak G\backslash \{0\}.$$ 
An iteration sequence is said to be {\it degenerate}  if there exists $n_0$ such that 
 $|z_n|=1$  for all $n\geq n_0$; it is said to be {\it nondegenerate}
 if it is not degenerate.

We note that for any $z \in \C'$, 
depending on the location of  $z$ in the unit sized square with integral
vertices that contains $z$, there are $2$, $3$ or all the $4$ of the vertices within
distance at most $1$ from $z$; $z_1$ can be chosen so that $|z_1|\geq 1$ and
$z_0-z_1^{-1}\in \frak G $ is one of these points, and  after the choice
of a $z_n$, $n\geq 1$, is made, between $2$ and $4$ choices are
possible for $z_{n+1}$, for continuing the iteration sequence.  For the points on the unit circle 
there are at least two nonzero 
points of $\frak G $ within distance $1$ that satisfy the nondegeneracy condition
at each $n$. Thus the sequence can be continued in
multiple ways at each stage to get a nondegenerate iteration sequence. 

\begin{remark}\label{algit}
{\rm Given an algorithm with the choice function $f:\C^* \to \frak G$ we get an iteration 
sequence $\{z_n\}$ by setting $z_0=z$ and $z_{n+1}=(z_{n}-f(z_n))^{-1}$ 
for all $n\geq 0$. The class of iteration sequences is however more general than
sequences arising from the algorithms.}
\end{remark}

\begin{remark}
{\rm Motivated by the classical simple continued
fractions,  one may consider stipulating that  $|z_{n}|>1$, for $n\geq 1$ but as
we shall see allowing $|z_{n}|\geq 1$ introduces no difficulty,
provided we ensure that $z_n-z_{n+1}^{-1}$ is nonzero, and this
enables somewhat greater generality. Also, since in $\C$
there exist irrational numbers $z$ with $|z|=1$ this generality seems
appropriate. 
}
\end{remark}

For each iteration sequence $\{z_n\}$ we get a sequence $\{a_n\}$ in
$\frak G$, called the corresponding {\it sequence of partial quotients}, defined by  
$$a_n=z_n-z_{n+1}^{-1}, \mbox { for all }n\geq 0; $$ 
we note that $a_n\neq 0$ for all $n \geq 1$. We shall show that if $z_n$
is nondegenerate and $\{p_n\}_{n=-1}^\infty$,
$\{q_n\}_{n=-1}^\infty$ is the $\cal Q$-pair corresponding to $\{a_n\}$
then $q_n\neq 0$ for all $n$ and $p_n/q_n$ converge
to $z$  (Proposition~\ref{prop2}). Thus corresponding to every
nondegenerate iteration sequence we get a continued fraction expansion for $z$. 
We shall also obtain an estimate for the difference $|z-\frac
{p_n}{q_n}|$ when the iteration sequence is bounded away from the unit
circle (see Proposition~\ref{prop2}).

\begin{proposition}\label{prop}
Let $z\in \C'$,  $\{z_n\}$ be an iteration sequence for $z$ and
$\{a_n\}$ be the corresponding sequence of partial quotients. Let  $\{p_n\}, \{q_n\}$
be the $\cal Q$-pair of sequences associated to $\{a_n\}$. Then 
 we have the following:  

i) $q_nz -p_n =(-1)^n(z_1\cdots z_{n+1})^{-1}$ for all $n\geq 0$,

ii) $
z=\displaystyle{\frac{z_{n+1}p_{n}+p_{n-1}}{z_{n+1}q_{n}+q_{n-1}}}$
for all $n\geq 0$,

iii) for all $n\geq 0$, $\displaystyle{|z-\frac {p_n}{q_n}|\leq \frac
1{|q_n|^2|z_{n+1}+\frac
      {q_{n-1}}{q_n}|}} $; in 
    particular, if $\{|q_n|\}$ is monotonically increasing and there
    exists a $\gamma >1$ such that 
    $|z_n|\geq \gamma $ for all $n$ then $\displaystyle{|z-\frac
      {p_n}{q_n}|\leq \frac 1{(\gamma-1)|q_n|^2}}$. 
\end{proposition}

\proof i) We shall argue by induction. Note that as $p_0=a_0$, $q_0=1$
and $z-a_0=z_1^{-1}$, the statement holds for $n=0$. Now let $n\geq 1$
and suppose by induction that the assertion holds for $0, 1, \dots, n-1$. 
Then we have 
$q_nz-p_n=(a_nq_{n-1}+q_{n-2})z-(a_np_{n-1}+p_{n-2})=a_n(q_{n-1}z-p_{n-1})+
(q_{n-2}z-p_{n-2})= (-1)^{n-1}(z_1\cdots 
z_n)^{-1}a_n+(-1)^{n-2}(z_1\cdots z_{n-1})^{-1}=(-1)^n(z_1\cdots 
z_n)^{-1}(-a_n+z_n)=(-1)^n(z_1\cdots z_{n+1})^{-1} $. This proves~(i). 

ii) By~(i) we have, for all $n\geq 1$, 
$$z_{n+1}= - \frac{q_{n-1}z  -p_{n-1}}{q_{n}z -p_{n}}.$$ 
Solving the equation for $z$ we get the desired assertion.

iii) Using  (iii) and the fact that $|p_{n-1}q_n-p_nq_{n-1}|=1$ we get 
 $$\displaystyle{|z-\frac
  {p_n}{q_n}|=|\frac{z_{n+1}p_{n}+p_{n-1}}{z_{n+1}q_{n}+q_{n-1}}-\frac
  {p_n}{q_n}|= \frac 1{|q_n|^2|z_{n+1}+\frac
      {q_{n-1}}{q_n}|}.} $$
The second assertion readily follows 
from the first.~\hfill $\Box$

\medskip
\begin{remark}\label{uncountable}
{\rm In the context of (iii) it may be noted here that we shall later
see sufficient conditions on $\{a_n\}$ which ensure that $\{|q_n|\}$ is
monotonically increasing. Also the condition $|z_n|\geq \gamma >1$
arises naturally under various circumstances, and holds in particular 
for iteration sequences arising from any algorithm whose fundamental set is contained in 
$B(1)$. It may be mentioned that the estimates we get from 
(iii) in Proposition~\ref{prop2} are rather  weak and not really significant 
from the point of view of Diophantine approximation; in fact even the pigeon hole principle
yields better inequalities. The main point about the statement  however is the
generality in terms of allowed sequences $\{a_n\}$ for which the
estimates hold. 

}
\end{remark}

\begin{remark}\label{degitsq} 
{\rm Let $R$
be the set of all $12$th roots of unity in $\C$. We note
that $R$ is precisely the set of points on the unit circle which are also at
unit distance from some nonzero Gaussian integer. (For $\rho =\pm 1$
and $\pm i$ there are three such points each, while for the other
$12\,$th roots there is a unique one.) Thus if $\{z_n\}$ is an iteration 
sequence for $z\in \C'$ and for some $m\in \N$ we have $|z_m|=|z_{m+1}|
=1$ then $z_m \in R$. It follows in particular that $\{z_n\}$ is a degenerate 
iteration sequence if and only if $z_n\in R$ for all large $n$. It may also be
seen, noting that $\{z_n\}$ is contained in $\C'$, that if, for some $m\in N$,
$z_m\in R$ and $|z_{m+1}| =1$ then $z_{m+1}=-z_m$.

}
\end{remark} 

\medskip
Towards proving that each nondegenerate iteration sequence of any $z\in \C'$ 
yields a continued fraction expansion of $z$ we first note the following:

\begin{proposition}\label{prop1}
Let $z\in \C'$ and $\{z_n\}$ be a nondegenerate iteration sequence for $z$. Then  
$\limsup |z_n|>1$.  
\end{proposition}

\proof Suppose this is not true. As $|z_n|\geq 1$ for all $n\geq 1$ this
means that $|z_n| \to 1$ as $n\to \infty$.  Then $|z_n-a_n|= |z_{n+1}|^{-1} 
\to 1$ as $n \to \infty $. As $a_n \in  \frak G \backslash \{0\}$ for all $n \geq 1$ 
the above two conditions imply that $d(z_n, R) \to 0$ as $n \to \infty$. 
Furthermore there exist $\rho_n \in R$ such that $|z_n-\rho_n| \to 0$. 
The conditions $|z_n-a_n| \to 1$ and $|z_n-\rho_n| \to 0$ imply in turn that
$|a_n-\rho_n| \to 1$ as $n \to \infty$.  Since $ \frak G$ and $R$ are discrete 
we get that $|a_n-\rho_n| = 1$ for all large $n$, say $n\geq n_0$.  
 We note also  that for any 
$\rho \in R$ and $a\in \frak G$ such that $|\rho -a|=1$, we have $\rho - a \in
R$; this 
can be seen from an inspection of the possibilities. In particular 
$\rho_n-a_n\in R$ for all  $n \geq n_0$. 

Now for all $n$ let $\zeta_n=z_n-\rho_n$ and $\beta_n=\rho_n-a_n$. Then we 
have $\beta_n\in R$ and  $z_n-a_n=\beta_n+\zeta_n$, and 
$\zeta_n\to 0$ as $n\to \infty$.
Then $$z_{n+1}=\frac 1{\beta_n+\zeta_n}=\frac 1{\beta_n} +\frac
{-\zeta_n}{\beta_n(\beta_n+\zeta_n)}.$$  Note that 
$$|\frac {\zeta_n}{\beta_n+\zeta_n}|\leq \frac{|\zeta_n|}{1-|\zeta_n|} \to 0.$$ 
Since $\beta_n\in R$, so is $\beta_n^{-1}$, and hence the two
relations above imply that, for all large $n$,  
$$\rho_{n+1}=\frac 1{\beta_n} \mbox{ and }\zeta_{n+1}= \frac
{-\zeta_n}{\beta_n(\beta_n+\zeta_n)}.$$ 
Then, for all large $n$, $$|\zeta_{n+1}|= |\frac
{\zeta_n}{\beta_n+\zeta_n}|=\frac {|\zeta_n|}{|z_n-a_n|} \geq |\zeta_n|,$$ 
since $|z_n-a_n|\leq 1$. 
Since   $\zeta_n\to 0$ as
$n\to \infty$ this implies that $\zeta_n= 0$ for all large $n$. 
Thus $z_n= \rho_n \in R$ for all large $n$. But this is a contradiction, since 
$\{z_n\}$ is a nondegenerate iteration sequence; see Remark \ref{degitsq}. 
Therefore the assertion as in the proposition must hold.~\hfill $\Box$ 

\medskip

\begin{theorem}\label{prop2}
Let $z\in \C'$,  $\{z_n\}$ be a nondegenerate iteration sequence for $z$ and
$\{a_n\}$ be the corresponding sequence of partial quotients. Let  $\{p_n\}, \{q_n\}$
be the $\cal Q$-pair of sequences associated to $\{a_n\}$. Let $n_0\geq 0$ 
be such that $|z_{n_0+1}|>1$. Then $q_n\neq 0$ for all $n\geq n_0$, and 
$\{p_n/q_n\}_{n_0}^\infty$ converges to $z$ as $n$ tends to infinity. In particular,
if $|z_1|$ or $|z_2|$ is greater than $1$ then $q_n\neq 0$ for all $n\geq 0$ 
and $\{a_n\}$ defines a continued fraction expansion for $z$.
\end{theorem}

\proof  By Proposition~\ref{prop}(i) we have $q_nz -p_n =(-1)^n(z_1\cdots z_{n+1})^{-1}$ 
for all $n\geq 0$. Under the condition in the hypothesis we get  that $0<|q_nz -p_n |<1$
for all $n\geq n_0$. Since $p_n\in \frak G$ this  shows in particular that $q_n\neq 0$ for 
all $n\geq n_0$. Furthermore, for $n\geq n_0$ we have $|z-\frac{p_n}{q_n}|= 
|z_1\cdots z_{n+1}|^{-1}|q_n|^{-1}\leq |z_1\cdots z_{n+1}|^{-1} \to 0,$
by Proposition~\ref{prop1}. The second part readily follows from the first,
noting that $q_0=1$. This proves the Theorem.  \hfill $\Box$

\section{Lagrange theorem for continued fractions}

In this section we prove an analogue of the classical Lagrange theorem
for the usual simple continued fraction expansions, that the expansion
is eventually periodic if and only if the number is a quadratic
surd.  We shall continue to follow the notation as before. 

A number $z\in \C$ is called a {\it quadratic surd} if it is not rational and 
satisfies a quadratic polynomial over $\frak G$, the ring of Gaussian
integers; namely, if there there exist $a,b,c \in \frak G$, with 
$a \neq 0$, such that 
$az^2+bz+c=0$; it may be noted that since $z$ is irrational such a polynomial
is necessarily irreducible over $\frak G$ and both its roots are irrational. 

Now, through the rest of the section let  $z\in \C'$ and
$\{z_n\}_{n=0}^\infty$ be an iteration  sequence for $z$. Also, as before let
$a_n=z_n-z_{n+1}^{-1} \in \frak G$,  for all $n\geq 0$, and $\{p_n\}$,
$\{q_n\}$ be the corresponding $\cal Q$-pair. 

\begin{proposition}\label{dir}
If $z_m=z_n$ for some  $0\leq m< n$, then $z$ is a quadratic
surd. 
\end{proposition}

\proof By Proposition~\ref{prop}(ii) we have
$$
z_m= -\frac{zq_{m-2}-p_{m-2}}{zq_{m-1}-p_{m-1}} \mbox{ and } 
z_n=-\frac{zq_{n-2}-p_{n-2}}{zq_{n-1}-p_{n-1}}. 
$$
It follows in particular that if $z_m$ is a quadratic surd then so is $z$. We
may therefore suppose that $m=0$, so $z_m=z$, and $n\geq 1$. The condition 
$z_n=z$ yields, from the above relation, that 
$z(zq_{n-1}-p_{n-1})+zq_{n-2}-p_{n-2}=0$, or $q_{n-1}z^2+(q_{n-2}-p_{n-1})z-p_{n-2}=0$. 
If $q_{n-1}\neq 0$ then this shows that $z$ is a quadratic surd and we are through. 
Now suppose that $q_{n-1}=0$. Since $q_0=1$, this means that $n\geq 2$. Also, 
as $q_{n-1}=0$, by Proposition~\ref{prop} we have $|p_{n-1}|=|z_1\cdots
z_n|^{-1}$. Since $p_{n-1}\in \frak G$ and $|z_k|\geq 1$ for all $k$, this yields 
that $|z_k|= 1$ for all $k=1, \dots , n$. Since $n\geq 2$, by Remark \ref{degitsq} 
it follows that $z_1\in R$. Furthermore, as it is irrational it must be one of the 
$12$th roots which is a quadratic surd. Hence $z$ is also a quadratic surd. 
 \hfill $\Box$

\medskip
We now prove the following partial converse of this.

\begin{theorem}\label{thm:lagr}
Let $z$ be a quadratic surd. Then for given $C>0$ there 
exists a finite subset $F$ of  $\C$ such that the following holds: if
$\{z_n\}$ is an  
iteration sequence for $z$, with the corresponding $\cal Q$-pair 
$\{p_n\}$, $\{q_n\}$, and 
$N$ is an infinite subset of $\N$ such that for all $n\in N$,
$$|q_{n-2}z-p_{n-2}|\leq C |q_{n-2}|^{-1}  \mbox{ and } |q_{n-1}z-p_{n-1}|\leq
C|q_{n-1}|^{-1},$$  
then there exists $n_0$ such that $\{z_n\mid n\in N, n\geq n_0\}$ is
contained in $F$.
\end{theorem}

\noindent{\it Proof}:
Let   $a,b,c \in\frak G$, with $a \neq 0$, be such that $ az^2+bz+c=0$. 
Let $C>0$ be given. Let $\cal P$ be the set of all quadratic
polynomials of the 
form $  \alpha z^2+\beta z+\gamma$, with $\alpha, \beta, \gamma \in
\frak G$ such that  $|\alpha|\leq C|2az+b|+1 , |\gamma |\leq
C|2az+b|+1 $ and $|\beta^2-4\alpha \gamma|= | b^2-4ac|$.  Then $\cal P$ is 
a finite collection of polynomials. Let $F$ be the set of all $\zeta
\in \C$ such that $P(\zeta)=0$ for some $P\in \cal P$; then $F$ is
finite. Now let $\{z_n\}$ be any   
iteration sequence for $z$ and $\{p_n\}$, $\{q_n\}$ be the
corresponding $\cal Q$-pair. By Proposition~\ref{prop2}(ii) the
condition implies that for all $n \geq 1$, 
$$
a(\frac{z_np_{n-1}+p_{n-2}}{z_nq_{n-1}+q_{n-2}})^2+
b\frac{z_np_{n-1}+p_{n-2}}{z_nq_{n-1}+q_{n-2}}+c=0\;\; \mbox{ or}
$$
$$
a(z_np_{n-1}+p_{n-2})^2+b(z_np_{n-1}+p_{n-2})(z_nq_{n-1}+q_{n-2})
+c(z_nq_{n-1}+q_{n-2})^2=0.
$$
We set
$$
A_n=ap_{n-1}^2+bp_{n-1}q_{n-1}+cq_{n-1}^2,\;\;C_n=A_{n-1} \;\;\mbox{ and}
$$
$$
B_n=2ap_{n-1}p_{n-2}+b(p_{n-1}q_{n-2}+q_{n-1}p_{n-2})+2cq_{n-1}q_{n-2}.
$$
Then we have 
$$
A_nz_n^2+B_nz_n+C_n=0 \mbox{ for all } n\geq 1.
$$
We note that  $A_n\neq 0$, since otherwise $p_{n-1}/q_{n-1}$ would be
a root of the equation 
$ax^2+bx+c=0$ which is not possible, since the roots of the polynomial are 
irrational.  We rewrite each
$A_n$ as 
{\begin{eqnarray}
\begin{array}{rll}
A_n &=& (ap_{n-1}^2+bp_{n-1}q_{n-1}+cq_{n-1}^2)-q_{n-1}^2(az^2+bz+c)\nonumber\\
  &=&a(p_{n-1}^2-q_{n-1}^2z^2) +b (p_{n-1}-q_{n-1}z)q_{n-1}\\
&=& (p_{n-1}-zq_{n-1})(a(p_{n-1}+zq_{n-1})+bq_{n-1})\\
&=& (p_{n-1}-zq_{n-1})(a(p_{n-1}-zq_{n-1})+(2az+b)q_{n-1}).\nonumber\\
\end{array} 
\end{eqnarray}}\\
As $|p_{n-1}-zq_{n-1}|=|z_1\cdots z_n|^{-1}\to 0$ as $n\to \infty$,
there exists $n_0$ such that for all $n\geq n_0$ in $N$ we have  
$$|A_n|\leq |2az+b||q_{n-1}||p_{n-1}-zq_{n-1}| +1\leq C|2az+b|+1,$$
by the condition in the hypothesis. Since $C_n=A_{n-1}$ we get also that 
$|C_n|\leq C|2az+b|+1$ for all large $n$ in $N$.
An easy  computation shows  that for all $n\geq 1$,
$B_n^2-4A_nC_n= b^2-4ac$. 
Thus for all $n\geq n_0$ in $N$, $A_nz^2+B_nz+C_n$ is a polynomial 
belonging to $\cal P$, and since $z_n$ is a root
of  $A_nz^2+B_nz+C_n$ we get that  $\{z_n\mid
n\in N, n\geq n_0\}$  is contained in $F$. This
proves the theorem. \hfill $\Box$ 

Theorem~\ref{thm:lagr} and Proposition~\ref{prop:cgcerate} imply the
following. 

\begin{corollary}\label{lagr-iter}
Let $z$ be a quadratic surd. Let  $\{z_n\}$ be an 
iteration sequence for $z$,  $\{a_n\}$ be the corresponding sequence,
in $\frak G$, of partial quotients,   and $\{p_n\}$, $\{q_n\}$ be the
corresponding 
$\cal Q$-pair. Suppose that  there exist $r_0\in \N$ 
and $\theta >1$ such that $|q_{n+r}|\geq \theta |q_{n}|$
for all $n\geq 0$ and $r\geq r_0$. Then $\{z_n\mid n\in \N\}$ is
finite. Consequently,   $\{a_n\mid n\in \N\}$ is also finite. 

\end{corollary}

Corollary~\ref{lagr-iter} may be viewed as the analogue Lagrange
theorem in terms of iteration sequences, under the mild growth
condition on the denominators of the convergents. Note that when
successive choices are not made according to any algorithm the 
sequence $\{a_n\}$ can not be expected to be eventually periodic. 
The growth condition on the denominators as above is satisfied 
is ensured under various conditions (see \S~6). 

\begin{remark}\label{rem1}
{\rm We have $
|q_nz
-p_n|=|q_n(\displaystyle{\frac{z_{n+1}p_n+p_{n-1}}{z_{n+1}q_n+q_{n-1}}-
\frac{p_n}{q_n}})|   
=\frac{1}{|z_{n+1}q_n+q_{n-1}|}$ for all $n\in \N$ and hence the condition
in the hypothesis of 
Theorem~\ref{thm:lagr} may be expressed as $|z_n+\frac
{q_{n-2}}{q_{n-1}}|>\delta \mbox{ and }|z_{n+1}+\frac 
{q_{n-1}}{q_{n}}|>\delta,$ for some $\delta >0$. In this respect it
may be observed that $|z_n|>1$ for all $n$ while $|q_{n-1}/q_n|$ may
typically be expected to be less than $1$; as $|p_n-q_nz|=|z_1\cdots
z_{n+1}|^{-1} \to 0$, $\{|q_n|\}$ is unbounded and hence the latter holds for 
infinitely many $n$, at any rate. The 
condition in the hypothesis would be satisfied if $\{|z_n|\mid n\in \N \}$ and
$\{|q_{n-1}/q_n|\mid n\in \N\}$ are at a positive distance from each other.
This motivates the following Corollary. 
}
\end{remark}

\begin{corollary}\label{cor:lagr}
Suppose that the $z$ as above is a quadratic surd. Let $N=\{n\in \N\mid
|q_{n-1}|>|q_{n-2}|\}$.  Suppose that there exists $\gamma>1$ such
that $|z_{n}|\geq\gamma$ for all $n\in N$. Then $\{z_n\mid n\in N\}$ is
finite.  
\end{corollary}

\proof For $n\in N$ we have  $|z_n+\frac
{q_{n-2}}{q_{n-1}}|\geq |z_n|-|\frac {q_{n-2}}{q_{n-1}}|>\gamma -1$. 
By the observation in 
Remark~\ref{rem1} we get $|q_{n-1}z-p_{n-1}|\leq C |q_{n-1}|^{-1}$,
with $C=(\gamma -1)^{-1}$.  
Thus if $n$ and $n+1$ are both contained in $N$, then 
 $|q_{n-1}z-p_{n-1}|\leq C |q_{n-1}|^{-1}$ and  $|q_nz-p_n|\leq
 C|q_{n}|^{-1}$. Now  suppose that $n\in N$ and $n+1\notin N$. Then
 $|q_n|\leq |q_{n-1}|$. Therefore $|q_nz-p_n|=|z_1\cdots z_{n+1}|^{-1}<
 |q_{n-1}z-p_{n-1}|\leq C 
 |q_{n-1}|^{-1}\leq C|q_{n}|^{-1}$. Thus the condition as in the
 hypothesis of Theorem~\ref{thm:lagr} is satisfied for all $n$ in $N$ as above. The corollary
 therefore follows from that theorem. \hfill $\Box$

We now apply the result to continued fraction expansions arising from
time-independent algorithms, say corresponding to a choice function
$f$. Let $z\in \C'$ and $\{z_n\}$ be the sequence produced from the
algorithm. We note that in this case if $z_n=z_m$ for some $m,n\in \N$
then $z_{n+l}=z_{m+l}$ for all $l\in \N$. 
We note also that if  $\Phi_f$ is contained in $B(0,r)$ for some $r<1$
then $|z_n|\geq r^{-1}$ for all $n\geq 1$. Let $\{a_n\}$ be the
corresponding continued fraction expansion. The sequence $\{a_n\}$ is
said to be {\it eventually periodic} if there exist $m, k\in \N$ such
that $a_{j+k}=a_j$ for all $j\geq m$. 

\begin{corollary}\label{lagrange}
Let $f$ be a choice function such that the fundamental set $\Phi_f$ is contained in 
$B(0,r)$ for some $r<1$. Then  $z\in \C'$ is a quadratic surd if and only if
its continued fraction expansion $\{a_n\}$ with respect to $f$ is eventually
periodic.  
\end{corollary} 

\proof Suppose $\{a_n\}$ is eventually periodic, and let $m, k\in \N$
be such that $a_{j+k}=a_j$ for all $j\geq m$. Then $z_{m+k}=z_m$, and
hence by Proposition~\ref{dir} $z$ is a quadratic surd. Now
suppose that $z$ is a quadratic surd. Then by
Corollary~\ref{cor:lagr} there exist  $0\leq  m< n$ such that 
$z_n=z_m$. Let $k=n-m$. Then $z_{m+k}=z_m$, and as noted above this
implies that $z_{j+k}=z_{j}$ for all $j\geq m$. Since $a_j=f(z_j)$ for
all $j\geq 0$ we get $a_{j+k}=a_{j}$ for all $j\geq m$. This shows that
$\{a_n\}$ is eventually periodic. \hfill $\Box$ 

\begin{remark}\label{nonperiodic}
{\rm Let $z\in \C'$ be a quadratic surd. There are at least two distinct eventually 
periodic continued fraction expansions for $z$ arising from two different 
algorithms, since given one such expansion we can find an algorithm to 
which the entries do not conform; furthermore the both the algorithms may
be chosen to be such that  their fundamental sets are contained in the open ball $B(0)$. 
Secondly, there are infinitely many eventually periodic expansions for $z$, 
since given any such expansions it can be altered after any given stage to get a new
one, which is also eventually periodic; not all of these may correspond to sequences
arising from algorithms. Clearly, the set of these eventually periodic expansions is
countable, since each of them is determined by two finite blocks
of Gaussian integers (one corresponding to an initial block and one to repeat itself in 
succession). 
Recall that for any $z$ there are uncountably many continued fraction expansions
(see Remark \ref{uncountable}).  Therefore there are also  uncountably many 
continued fraction expansions which are not eventually periodic.  
}
\end{remark}

\begin{remark}
{\rm We note that the eventually periodic continued fraction expansions
of a quadratic surd corresponding to two different algorithms may differ in their 
eventual periodicity as also in respect of 
whether the expansion is purely periodic (periodic from the beginning) 
or not. For example for $\frac{\sqrt 5+1}2$ we have the usual expansion 
$\{a_n\}$ with $a_n=1$ for all $n$, which may  also be seen to be its continued fraction expansion 
with respect any of the shifted Hurwitz algorithms $f_d$ (see Example 2 in \S 2), for
any real number $d$ in the interval $(\frac{\sqrt 5-1}2-\frac12, 1)$; we recall
that from among these for $d\in (\frac{\sqrt 5-1}2-\frac12, 1-\frac 1{\sqrt 2})$
the associated set 
$\Phi_f$ is contained in $B(0,r)$ for some $r<1$, thus satisfying the condition in 
Corollary \ref{lagrange}. On the other hand, for the continued 
fraction expansion of $\frac{\sqrt 5+1}2$ with respect to the Hurwitz algorithm
may be seen to be given by $\{a_n\}$ where $a_0=2$ and 
$a_n=(-1)^n3$ for all $n\geq 1$; it is not purely periodic, and the eventual
period is 2, rather than 1. 
}
\end{remark}

\section{Growth properties of the convergents}

In this section we discuss various growth properties of the sequences $\{q_n\}$
from the $\cal Q$ pairs. We begin with the following general fact.

\begin{proposition}\label{alternate}
Let $\{a_n\}_{n=0}^m$, where $m\in \N$, be a sequence in $\frak G$ such that 
$|a_n|>1$ for all $n\geq 1$
and let $\{p_n\}$, $\{q_n\}$ be the corresponding $\cal Q$-pair. Suppose that 
$q_n\neq 0$ for all $n\geq 0$. Let $\sigma \in 
(\sqrt 2, \frac {\sqrt 5 +1}2)$, and   $1\leq n\leq m-1$
be such that $|{q_{n-1}}|> \sigma |{q_{n-2}}|$. Then the following statements hold: 
i) if  $|a_n|>2$ then $|{q_{n}}|> \sigma |{q_{n-1}}|$; ii) if  $|a_n|\leq 2$ and 
$\hbox{\rm Re } a_na_{n+1}\geq \chi$, where $\chi =2 $ if $|a_{n+1}|=\sqrt 2$ and $0$
otherwise, then  $|{q_{n+1}}|> \sigma|{q_{n}}|$. 


\end{proposition}

\proof 
Let $r_k=q_{k-1}/q_k$, for all $k\geq 0$,  and $r=\sigma^{-1}$. Then 
$r\in (\frac {\sqrt 5-1}2, 1/\sqrt 2)$. 
We have $r_n=1/(a_n+r_{n-1})$, so $|a_n|\leq |r_{n-1}|+|r_{n}|^{-1}<
r+|r_{n}|^{-1}$. When $|a_n|>2$  we get  
$r+|r_{n}|^{-1}> |a_n|\geq \sqrt 5> r+r^{-1}$,  since $r>\frac {\sqrt
  5-1}2$ and $r+r^{-1}$ is a decreasing function in $(0,1)$. Therefore   
$|r_n|<r$. This proves 
the first statement. Now  suppose that $|a_n|\leq 2$ and 
 $\hbox{\rm Re } a_na_{n+1}\geq \chi$ where
$\chi $ is as above; in particular $\hbox{\rm Re } a_na_{n+1}\geq 0$.

Before continuing with the proof we note that for any $r>0$ and $z\in \C$  such that 
$|z|>r$ the set $B(z,r)^{-1}:=\{\zeta^{-1}\mid \zeta \in B(z,r)\}$ of the inverses of 
elements from $B(z,r)$, 
is precisely $\displaystyle B(\frac {\bar z}{|z|^2-r^2}, \frac {r}{|z|^2-r^2})$. 

Now, as $r_{n-1}\in B(0,r)$ and $r<|a_n|$, we get that $r_n$, which is $1/(a_n+r_{n-1})$,
is contained in $B(a_n,r)^{-1}= B(y{\bar a_n}, {yr})$, where $y=1/(|a_n|^2-r^2)$. 
We have $r<1/\sqrt 2$, so $ 1/(\sqrt 2-r)<1/(\sqrt 2-1/\sqrt 2)=
\sqrt 2$. Since $|a_n|$ and $|a_{n+1}|$ 
are at least $\sqrt 2$, this shows that  $|a_{n+1}|>1/({|a_n|-r})$, 
which is the same as ${y(|a_n|+r)}$. Therefore $|y\bar a_n + a_{n+1}|\geq y|a_n|-|a_{n+1}|
>yr$. We have 
$r_{n+1}=  1/ (a_{n+1}+r_n)\in B(a_{n+1}+y\bar a_n, yr)^{-1} $ and in the light of  the 
observation above it follows that $r_{n+1}\in B(\frac {\bar z}{|z|^2-s^2}, \frac {s}{|z|^2-s^2})$, 
where $z=y\bar a_n +a_{n+1}$ and 
$s=yr$. We would like to conclude that $r_{n+1}\in B(0, r)$ and 
for this it suffices to show that $|\frac{\bar z}{|z|^2-s^2}|+ |\frac {s}{|z|^2-s^2}|<r$, or 
equivalently that $r(|z|-s)>1$. Substituting for $z$ and $s$ and putting $x=y^{-1}$
we see that the above condition is equivalent to $r(|\bar a_n+xa_{n+1}|-r)>x=|a_n|^2-r^2$, 
or equivalently $r|\bar a_n +xa_{n+1}|>|a_n|^2$. Squaring both the sides and recalling
that  $r^2=|a_n|^2-x$ the condition may be written as 
$$(|a_n|^2-x)(\bar a_n+xa_{n+1}) (a_n+x\bar a_{n+1})>|a_n|^4.$$
The condition can be readily simplified to 
$$|a_{n+1}|^2x^2+(-\alpha +\beta)x-|a_{n}|^2(\beta-1)<0,$$ 
where $\alpha = |a_{n}|^2|a_{n+1}|^2$ and 
$\beta=2\hbox{\rm Re } a_na_{n+1}$. The discriminant $\Delta$ of the quadratic polynomial 
is seen to be $(\alpha +\beta)^2-4\alpha$. Since $\beta =2\hbox{\rm Re } a_na_{n+1}\geq 0$ 
and $\alpha \geq 4$, we get that $\Delta \geq 0$.

 The   (real) roots of the quadratic are given by 
$\lambda_\pm =(\alpha -\beta \pm \sqrt \Delta)|a_{n}|^{2}/2\alpha$ and the desired 
condition as above 
holds if $x\in (\lambda_-, \lambda_+)$. As $r^2=|a_{n}|^{2}-x$ we see that 
the condition holds if $r^2$ is in the interval $(|a_{n}|^{2}-\lambda_+, 
|a_{n}|^{2}-\lambda_-)$, or equivalently if the following inequalities hold:
$$(2\alpha)^{-1}|a_{n}|^{2}(\alpha +\beta-\sqrt \Delta)
<r^2 < (2\alpha)^{-1}|a_{n}|^{2}(\alpha +\beta+\sqrt \Delta).$$ We note that the last term 
is at least $(2\alpha)^{-1}|a_{n}|^{2}\alpha= |a_{n}|^{2}/2\geq 1$, and hence  it is greater 
than $r^2$. It now suffices to prove that the first inequality holds, or
equivalently that
$$\sigma^2 < (2\alpha)|a_{n}|^{-2}(\alpha +\beta-\sqrt \Delta)^{-1}= \frac 12 |a_{n}|^{-2}
(\alpha +\beta+\sqrt \Delta),$$
as $(\alpha +\beta)^2-\Delta =4\alpha$. 
Now first suppose that $|a_{n+1}|\geq 2$. Then, as $\beta \geq 0$,
the right hand side is at least $\frac 12 |a_{n}|^{-2} (4  |a_{n}|^2
+ \sqrt{16 |a_{n}|^4 -16 |a_{n}|^2})= 2(1+\sqrt {1-|a_n|^{-2}})\geq 2+\sqrt 2$, and hence 
it indeed exceeds $\sigma^2$.
Now suppose $|a_{n+1}|=\sqrt  2$. Then we have $\hbox{\rm Re } a_na_{n+1}\geq 2$. 
Under these conditions it can  be explicitly verified that the right hand expression above 
exceeds $(\frac {\sqrt 5+1}2)^2=\frac {3+\sqrt 5}2$; for instance,  if $a_{n+1}=1+i$ then 
the only possibilities for $a_n$ in $\frak G$ with
$|a_n|\leq 2$ and   $\hbox{\rm Re } a_na_{n+1}\geq 2$ are $1-i, 2$ or $-2i$ and the 
values of the expression may be seen to be $2+\sqrt 3$ in the first case and $(3+\sqrt 7)/2$
in the other two cases; the corresponding assertion  for other choices
of $a_{n+1}$ 
follows from this from symmetry considerations, or may be verified
directly.
Therefore $r^2$ is contained in the desired interval and
the argument shows that $r_{n+1}\in B(0,r)$. This proves the proposition. \hfill $\Box$

\medskip
The following corollary follows immediately, by application of Proposition~\ref{alternate}
to values of $\sigma$ smaller than $(\sqrt 5+1)/2$. 

\begin{corollary}\label{cor:alternate}
Let $\{a_n\}_{n=0}^m$ be a sequence in $\frak G$ such that $|a_n|>1$ for all $n\geq 1$
and let $\{p_n\}$, $\{q_n\}$ be the corresponding $\cal Q$-pair. Suppose that 
$q_n\neq 0$ for all $n\geq 0$. Let  $1\leq n\leq m-1$
be such that ${|{q_{n-1}}|\geq (\frac {\sqrt 5 +1}2) |{q_{n-2}}|}$, and either $|a_n|>2$ or 
$\hbox{\rm Re } a_na_{n+1}\geq \chi$, where $\chi =2 $ if $|a_{n+1}|=\sqrt 2$ and $0$
otherwise. Then 
either $|{q_{n}}|/|{q_{n-1}}|$ or $|{q_{n+1}}|/|{q_{n}}|$ is at least $\frac {\sqrt 5+1}2$. 

\end{corollary}

If $\{a_n\}$ is a sequence in $\frak G$ corresponding to continued fraction expansion 
with respect to the Hurwitz algorithm then for all $n\geq 1$ such that $|a_n|\leq 2$ we have
$\hbox{\rm Re } a_na_{n+1}\geq \chi$, where $\chi$ is as in Proposition~\ref{alternate}
(see \cite{H} or \cite{Hen}). We thus get the following Corollary. 

\begin{corollary}\label{Hur:alternate}
Let $\{a_n\}$ be a sequence in $\frak G$ given by the continued fraction expansion of
a complex number with respect to the Hurwitz algorithm, and let  $\{p_n\}$, $\{q_n\}$ be the corresponding $\cal Q$-pair.  Then for any $n\geq 1$ either $|q_n|/|{q_{n-1}}|$ or 
$|{q_{n+1}}|/|{q_{n}}|$ is at least ${(\sqrt 5+1)}/2$. 
\end{corollary}

\proof We note that the condition in Corollary \ref{cor:alternate} holds for $n=1$, since
$q_{-1}=0$. Applying the Corollary successively we get that for any $n\geq 1$ one of the 
ratios $|q_n|/|{q_{n-1}}|$ or 
$|{q_{n+1}}|/|{q_{n}}|$ is at least ${(\sqrt 5+1)}/2$. \hfill $\Box$
 
 \begin{remark}{\rm
 Corollary \ref{Hur:alternate} strengthens  Lemma~5.1 of \cite{Hen}, where 
 the corresponding assertion is proved (for sequences arising from the Hurwitz algorithm) 
 with the constant $\frac 32$ in place of  $\frac{\sqrt 5+1}2$.
 }
 \end{remark}
 
As noted earlier, in various contexts it is useful to know whether 
$|q_n| $ is increasing with $n$. We shall next address the issue and describe conditions
under which this is assured.  In this respect we first note the following simple fact about
certain finite sequences of complex numbers, which will also be used later.  

We denote by $\sigma_y:\C \to \C$ the reflexion in the $y$-axis, namely the map defined 
by $\sigma_y(z)=-\bar z$ for all $z \in \C$.

\begin{lemma}\label{simple}
Let  $z_0, z_1, \dots , z_{m}$ be nonzero complex numbers such that 
 $|z_0|\leq 1$,  $|z_n| > 1$ for $n=1, \dots , m$, and let $b_n=z_n-z_{n+1}^{-1}$ 
 for all $n=0, \dots , m-1$.  Suppose that $b_n\in \frak G$ for all
 $n\leq m-1$,  $|b_0| >1$ and that for  $1\leq n\leq m-1$ either 
 $|b_n-\sigma^n_y(b_{0})| \geq  2$ or   
 $b_n= 2\sigma^n_y(b_{0})$.
Then $|b_0|=\sqrt 2$, $b_n= 2\sigma_y^n(b_0)$  for $n=1, \dots , m-1$,  and  
 $|z_{n+1}|\leq \sqrt 2 +1$  for all $n=1, \dots , m$.
  
\end{lemma}

\proof We have $b_0=z_0-z_{1}^{-1}$, $|z_0|\leq 1$ and $|z_1|> 1$, which 
yields $|b_0| < 2$, and since $b_0 \in \frak G$ and $|b_0|>1$ we get $|b_0|=  \sqrt 2$.
 We have 
$z_1^{-1}\in \bar B(-b_0) \cap B$ and hence we get that  $z_1 \in 
\bar B(\sigma_y(b_0)) \backslash \bar B$, and in particular $|z_1|\leq \sqrt 2 +1$. Also, 
since $z_1 \in \bar B(\sigma_y(b_0))$ and $|z_1-b_1|<1$ we get 
that $|b_1-\sigma_y(b_0)| <2$.  Hence by the condition in the
hypothesis  $b_1= 2\sigma_y(b_{0})$. 
Then $z_2^{-1}+2\sigma_y(b_0)=z_1\in \bar B(\sigma_y(b_0))$, 
and so $z_2^{-1}\in \bar B (-\sigma_y(b_0))$. Also, by hypothesis we have
$z_2^{-1} \in B$. Hence $z_2\in 
\bar B (\sigma^2_y(b_0))\backslash \bar B$, and in particular $|z_2|\leq \sqrt 2 +1$.
Again, since $z_2 \in \bar B(\sigma_y^2(b_0))$ and $|z_2-b_2|<1$ this
implies  
that $|b_2-\sigma^2_y(b_0)| <2$. Therefore by the condition in the
hypothesis we get that  $b_2= 2\sigma^2_y(b_{0})$.  
Repeating the argument we see $b_n= 2\sigma_y^n(b_0)$   and
$|z_{n+1}|\leq \sqrt 2 +1$  for all $n=0, \dots, m-1$. This proves the
Lemma.\hfill $\Box$ 

\medskip
\begin{definition}{\rm  A sequence $\{a_n\}_{n=0}^\infty$ in $\frak G$
    such that $|a_n|>1$  for all 
$n\geq 1$ is said to satisfy {\it Condition (H)} if for all $i,j$
such  that $1\leq i <j$,  
$|a_j|=\sqrt 2$ and $a_{k}=2\sigma_y^{j-k}( a_j)$ 
for all $k\in \{i+1,\dots, j-1\}$ (empty set if $i=j-1$), either
$a_{i}=2\sigma_y^{j-i}( a_j)$  or $|a_{i}-\sigma_y^{j-i} (a_j)|\geq  2$ .
}
\end{definition}

 (It may be noted that depending on the specific value of 
$a_j$ and whether $j-i$ is odd or even, in the above definition, $a_i$ is disallowed to 
belong to a set of five elements of $\frak G$ with absolute
value greater than $1$; if $a_j=1+i$ and $j-i$ is odd then these are $2i, -1+i, -1+2i, -2$ 
and $-2+i$; 
for other $a_j$ with  $|a_j|=\sqrt 2$ the corresponding sets can be written down using
symmetry under the maps $z\mapsto -z$ and $z\mapsto \bar z$).

\begin{definition}{\rm  A sequence $\{a_n\}_{n=0}^\infty$ in $\frak G$ such that $|a_n|>1$  for all
$n\geq 1$ is said to satisfy {\it Condition (H$'$)} if  $|a_{n}+\bar a_{n+1}|\geq 2$ for all $n\geq 1$ 
such that $|a_{n+1}|=\sqrt 2$.
}
\end{definition}

\begin{remark}{\rm 
We note that if  Condition (H$'$)  holds, then Condition (H) also holds, trivially. 
Condition (H) is strictly weaker than Condition (H$'$); viz. 
it can be seen there exists $z\in \C'$ such that in the corresponding
continued fraction sequence $\{a_n\}$ via the nearest integer
algorithm $-2+2i$ may be followed by $1+i$. While Condition (H$'$) is 
simpler, and may be used with other algorithms (see Theorem~\ref{represent}), 
Condition (H)   is of significance 
since the sequences occurring as continued fraction expansions with respect
the Hurwitz algorithm satisfy it, 
as was observed in Hurwitz \cite{H}, but not Condition (H$'$) in general.  
}
\end{remark}

\begin{proposition}\label{monotone}
Let $\{a_n\}_{n=0}^\infty$ be a sequence in $\frak G$, with $|a_n|>1$  for all
$n\geq 1$, satisfying Condition (H) and  let $\{p_n\}$, $\{q_n\}$ be the corresponding 
$\cal Q$-pair.  Then $|q_n|> |q_{n-1}| $  for all $n\geq 1$. 
\end{proposition}

\proof If possible suppose $|q_{n+1}/q_{n}|\leq 1$ for some $n\in \N$ and let $m$ 
be the smallest integer for which this holds. For all $0\leq k\leq m$ let 
$z_k=q_{m-k+1}/q_{m-k}$, and $b_k=a_{m-k+1}$. We have $|z_0|\leq 1$, while  
 for all $k=1, \dots, m$, $|z_k|> 1$. Also $z_k=\displaystyle{\frac{q_{m-k+1}}{q_{m-k}}
=a_{m-k+1}+\frac{q_{m-k}}{q_{m-k-1}}= a_{m-k+1}+z_{k+1}^{-1}=b_{k}+z_{k+1}^{-1}}$
for all $k\leq m$. As $\{a_n\}$ satisfies Condition~(H) it follows
that the condition in  Lemma~\ref{simple} is satisfied for $z_0, z_1,
\dots , z_m$ as above, and the lemma therefore implies 
that $|b_0|=\sqrt 2$, $b_m= 2\sigma_y^m(b_0)$ 
and $|z_m|\leq \sqrt 2 +1$. But this is a contradiction, since 
$z_m=q_1/q_0=a_1=b_m$ and $|b_m|=|2\sigma_y^m(b_0)|=2\sqrt 2$ . Therefore
 $|q_n|>|q_{n-1}|$ for all $n\geq 1$. This proves the proposition. \hfill $\Box$

\begin{definition}
{\rm A sequence $\{a_n\}$ in $\frak G$ such that $|a_n|>1$ for all
$n\geq 1$, is said to satisfy {\it Condition (A)} if for all $n\geq 1$ such that   $|a_n|\leq 2$,  
$\hbox{\rm Re } a_na_{n+1}\geq \chi$, where $\chi =2 $ if $|a_{n+1}|=\sqrt 2$ and $0$
otherwise. (If $|a_n| \leq 2$ then $a_{n+1}$ is required to be such that $\hbox{\rm Re } 
a_na_{n+1}\geq 0$; if moreover $|a_{n+1}|=\sqrt 2$ then the condition further stipulates
that $\hbox{\rm Re } a_na_{n+1}\geq 2$).
}
\end{definition}

We  note that both  Condition (H) and Condition (A)  hold 
for the sequences $\{a_n\}$ arising from the Hurwitz algorithm (see \cite{H}, \cite{Hen}). 
Corollary \ref{cor:alternate} and Proposition~\ref{monotone} together imply the following. 

\begin{proposition}\label{growth}
Let $\{a_n\}$ be a sequence in $\frak G$, with $|a_n|>1$ for all
$n\geq 1$, satisfying Conditions (H) and (A), and let $\{p_n\}$, $\{q_n\}$ be 
the corresponding $\cal Q$-pair.  Then $|q_n|> (\frac {\sqrt 5 +1}2) |q_{n-2}| $ for all 
$n\geq 2$. 
\end{proposition}

An analogous conclusion with a smaller constant can also be proved under
a weaker condition than Condition (A), which has also the advantage that it  
involves exclusion of occurrence of only finitely 
many pairs, as $(a_n, a_{n+1})$.

\begin{definition}{\rm
A sequence $\{a_n\}$ in $\frak G$ such that $|a_n|>1$ for all
$n\geq 1$, is said to satisfy {\it Condition (A$'$)} if for all $n\geq 1$ such that $|a_{n}|=\sqrt 2$ 
and $|a_{n+1}|\leq  2$, 
$a_{n+1}$ is one of the Gaussian integers $\bar a_n -a_n, \bar a_n $ or $\bar a_n +a_n$.
}
\end{definition}

(The condition can also be stated equivalently as follows: if $|a_{n}|=\sqrt 2$ 
and $|a_{n+1}|\leq  2$, then $|a_na_{n+1}+1|\geq |a_{n+1}|+1$; this  may be explicitly 
verified to be equivalent to the above condition; 
when the latter inequality holds, in fact one has $|a_na_{n+1}+1|> |a_{n+1}|+\sqrt 5 -1$.)

\begin{proposition}\label{growth2}
Let $\{a_n\}$ be a sequence in $\frak G$ with  $|a_n|>1$ for all
$n\geq 1$, satisfying Conditions (H) and (A$'$), and let $\{p_n\}$, $\{q_n\}$ 
be the corresponding $\cal Q$-pair. 
  Then $|q_n|> (\sqrt 5-1) |q_{n-2}| $ for all $n\geq 1$. 
\end{proposition}

 \proof Since $q_{-1}=0$ and $q_1=a_1$ the assertion trivially holds for $n=1$.
 Now suppose that $n\geq 2$. Since $\{a_n\}$ satisfies Condition (H), by 
 Proposition~\ref{monotone} $\{|q_n|\}$ is increasing, and in particular $q_n\neq 0$
 for all $n$. Recall that, we have $q_n=a_nq_{n-1}+q_{n-2}$ and
hence $|q_n/q_{n-1}|\geq |a_n|-|q_{n-2}/q_{n-1}|>|a_n|-1$, as $|q_{n-2}/q_{n-1}|<1$. 
Since $|q_{n-1}|>|q_{n-2}|$, when 
$|a_n|\geq \sqrt 5$ this yields $|q_n/q_{n-2}|>|q_n/q_{n-1}|>\sqrt 5 -1$, as
desired. Similarly if 
$|a_{n-1}|\geq \sqrt 5$ then we have
$|q_n/q_{n-2}|>|q_{n-1}/q_{n-2}|>\sqrt 5 -1$.  Now suppose that both 
$|a_n|$ and $|a_{n-1}|$  are less than $\sqrt 5$, namely
$|a_{n-1}|\leq 2$ and $|a_n|\leq 2$. We recall that
$q_n=a_nq_{n-1}+q_{n-2}=a_n(a_{n-1}q_{n-2}+q_{n-3}) 
+q_{n-2}=(a_na_{n-1}+1)q_{n-2}+a_nq_{n-3}$,  
and hence $|q_n/q_{n-2}|>|a_na_{n-1}+1|-|a_n|$. When $|a_{n-1}|=2$ it 
can be explicitly verified that for all possible choices of  $a_n$ with $|a_n|\leq 2$
that are admissible under Condition (H) we have $|a_na_{n-1}+1|-|a_n|>\sqrt 5 -1$. 
Now suppose $|a_{n-1}|=\sqrt 2$. Then as  noted above,  
Condition (A$'$) implies that for all admissible choices of $a_n$, 
$|a_na_{n-1}+1|-|a_n|$ exceeds $\sqrt 5 -1$.  This proves the proposition. \hfill $\Box$ 

\begin{remark}
{\rm If $\{a_n\}$ is a sequence with $a_n\in \Z \subset \frak G$ for all $n\geq 0$ and 
$|a_n|\geq 2$ for all $n\geq 1$, then Condition (H) and Condition (A$'$) are satisfied
and hence the conclusion of Proposition~\ref{growth2} holds for these sequences. If 
moreover  $a_{n+1}/a_n$ is positive whenever $|a_n|=2$, then Condition (A) also 
holds and in this case the conclusion of Proposition~\ref{growth} holds. 
}
\end{remark}

\section{Sequence spaces and continued fraction expansion}

In this section we introduce various spaces of sequences which serve as continued 
fraction expansions. 

 \medskip
For $\theta>1$ let 
$\Omega_E^{(\theta)}$ denote the set of sequences $\{a_n\}$ such that if $\{p_n\}$,
$\{q_n\}$ is the corresponding $\cal Q$-pair then,  for all $n\geq 1$, $|q_n|> |q_{n-1}| $ 
and $|q_n|\geq \theta |q_{n-2}| $; we note that for these sequences $|q_n|\geq \theta^{(n-1)/2} $
for all $n\geq 0$ (the suffix $E$  in the notation signifies exponential growth of the $|q_n|$'s). 
By Proposition~\ref{growth} the sequences  $\{a_n\}$ satisfying Condition (H) and
Condition (A) are contained in $\Omega_E^{(\frac {\sqrt 5+1}2 )}$; in particular every 
sequence arising as the continued fraction expansion of any $z$ in $\C'$ with respect 
to the Hurwitz algorithm is contained in $\Omega_E^{(\frac {\sqrt 5+1}2 )}$. Similarly 
by Proposition~\ref{growth2} the sequences  $\{a_n\}$ satisfying Condition (H) and
Condition (A$'$) are contained in $\Omega_E^{(\sqrt 5-1)}$. 

\begin{theorem}\label{cgce}
Let  $\{a_n\}_{n=0}^\infty$ be  a sequence in $\Omega_E^{(\theta)}$, for some $\theta >1$,
 and let  $\{p_n\}, \{q_n\}$ 
the corresponding $\cal Q$-pair.  Then the following conditions hold:

i) $\{p_n/q_n\}$  converges as $n\to \infty$;

ii) if $z$ is the limit of $\{p_n/q_n\}$ then ${\displaystyle |z- \frac{p_n}{q_n}|\leq \frac c{|q_n|^2}}$, 
where ${\displaystyle c=\frac{2\theta^2}{(\theta^2 -1)}}$;

iii) the limit $z$ is an irrational complex number, viz. $z\in \C'$.
\end{theorem} 

\proof We have $p_{n+1}q_{n}-p_{n}q_{n+1}=(-1)^{n}$ for all
$n\geq 0$, so
$\displaystyle \frac {p_{n+1}}{q_{n+1}}-\frac {p_{n}}{q_{n}} =\frac
  {(-1)^{n}}{q_{n}q_{n+1}}$, and in turn 
$\displaystyle |\frac {p_{n+m}}{q_{n+m}}-\frac {p_{n}}{q_{n}}|\leq
\sum_{k=0}^{m-1} \frac   {1}{|q_{n+k}q_{n+k+1}|}$ for all $n\geq 0$ and 
$m\geq 1$. As $|q_k/q_{k-2}|\geq \theta$ for all $k\geq 1$  this shows that
$\{p_n/q_n\}$ is convergent.  This proves~(i). Assertion~(ii) then follows 
from Proposition~\ref{prop:cgcerate}. Finally, it is easy to see
that such an approximation can  
not hold for rational complex numbers and hence the limit must be 
an irrational number. \hfill $\Box$

\medskip
 For any $\{a_n\}$  in $ \Omega_E^{(\theta)}$, $\theta>1$, 
we denote by $\omega (\{a_n\})$ the limit 
of  $\{p_n/q_n\}$ where $\{p_n\}, \{q_n\}$ is the corresponding $\cal Q$-pair. 
The following corollary shows that for any fixed $\theta>1$ the convergents 
coming from possibly different representations of a number $z$ from $\Omega_E^{(\theta)}$
also converge to $z$, provided the corresponding indices go to infinity.

\begin{corollary}\label{commoncf}
Let $\theta >1$, $z\in \C'$ be given, and  
$\{l_j\}$ be a sequence in $\N$ such that $l_j\to \infty $ as $j\to \infty$. Let $\{z_j\}$ be 
a sequence in $\C$ such that the following holds: for each $j$ there exists a
sequence $\{a_n\}$ in  $\Omega_E^{(\theta)}$ such that $z=\omega (\{a_n\})$ and $z_j=
p_n/q_n$ for some $n\geq l_j$, where $\{p_n\}$, $\{q_n\}$ is the $\cal Q$-pair corresponding 
to $\{a_n\}$. Then $z_j\to z$ as $j\to \infty$.
\end{corollary}

\proof By Proposition~\ref{cgce} we have ${\displaystyle |z-z_j|\leq 
\frac c{|q_n|^2}\leq \frac c{\theta^{n-1}} \leq \frac c{\theta^{l_j-1}}}$, where 
$c={{2\theta^2}/{(\theta^2 -1)}}$. As  $l_j\to \infty$ this shows 
that $z_j\to z$, as $j\to \infty$. \hfill $\Box$

\medskip
We equip $\Omega_E^{(\theta)}$, $\theta>1$, with the topology induced as a subset of the space
of all sequences $\{a_n\}$ with $a_n\in \frak G$ and $|a_n|>1$ for all $n\geq 1$,
the latter being equipped with the product topology as a product of the individual 
copies of $\frak G$ (for the $0$'th component) and $\frak G\backslash \{0\}$ (for 
the later components) with the discrete topology. 

\begin{proposition}\label{cont}
For all $\theta >1$ the map $\omega : \Omega_E^{(\theta)}\to \C$ is continuous. 
\end{proposition}

\proof 
Let  $\{a_n\}$  be a sequence from $\Omega_E^{(\theta)}$  and let $\{p_n\}, \{q_n\}$ be the 
corresponding $\cal Q$-pair. Let $\epsilon >0$ be given and $m\in \N$ be such 
that $|q_m|^2> {2c}\epsilon^{-1}$, where $c={{2\theta^2}/{(\theta^2 -1)}}$.
If  $\{a_n'\}$  is any sequence in $\Omega_E^{(\theta)}$
such that $a_n'=a_n$ for $n=0, 1, \dots , m$, and $\{p_n'\}, \{q_n'\}$ is the 
corresponding $\cal Q$-pair then $p_n'=p_n$ and $q_n'=q_n$ for $n=0, 1, \dots , m$. 
Then by Proposition~\ref{cgce} (ii), $\displaystyle{| \omega (\{a_n\}) - \frac{p_m}{q_m}|\leq 
\frac{c}{|q_m|^2}}< \epsilon/2$, as well as $\displaystyle{| \omega (\{a_n'\}) - 
\frac{p_m}{q_m}|}<
\epsilon/2$, so $| \omega (\{a_n'\})- \omega (\{a_n\})| < \epsilon$. This shows that
$\omega$ is continuous. \hfill $\Box$

\bigskip
The spaces $\Omega_E^{(\theta)}$ deal with continued fraction expansions in terms 
of the growth of the $|q_n|$'s, and do not seem to be amenable to a priori description, 
though of course as seen earlier verification of certain conditions on the sequences 
enables to conclude that they belong to  $\Omega_E^{(\theta)}$ for certain values of $\theta$.
Our next objective is to describe a space of sequences defined in terms of nonoccurrence of 
certain simple blocks, related to Conditions (H$'$) and (A$'$), which yield continued 
fraction expansion. We begin by introducing an algorithm such that the corresponding
sequences of partial quotients satisfy these conditions. 

Let $S$ be the square $\{z=x+iy\mid |x|\leq \frac 12,|y|\leq \frac 12  \}$. Let $P$ be the 
set of $8$ points on the boundary of $S$ such that one of the coordinates is $\pm \frac 16$
(the other being $\pm \frac 12$). Let $R$ be the set defined by
$$R=\{z\in S\mid d(z,p)\geq \frac 13 \hbox{ \rm for all } p \in P\}.$$
It may be noted that $R$ is a star-like subset with the origin as the
center, bounded by $8$ circular segments; see Figure~1.   

\begin{figure}[h]
\centering

\includegraphics[scale=0.8]{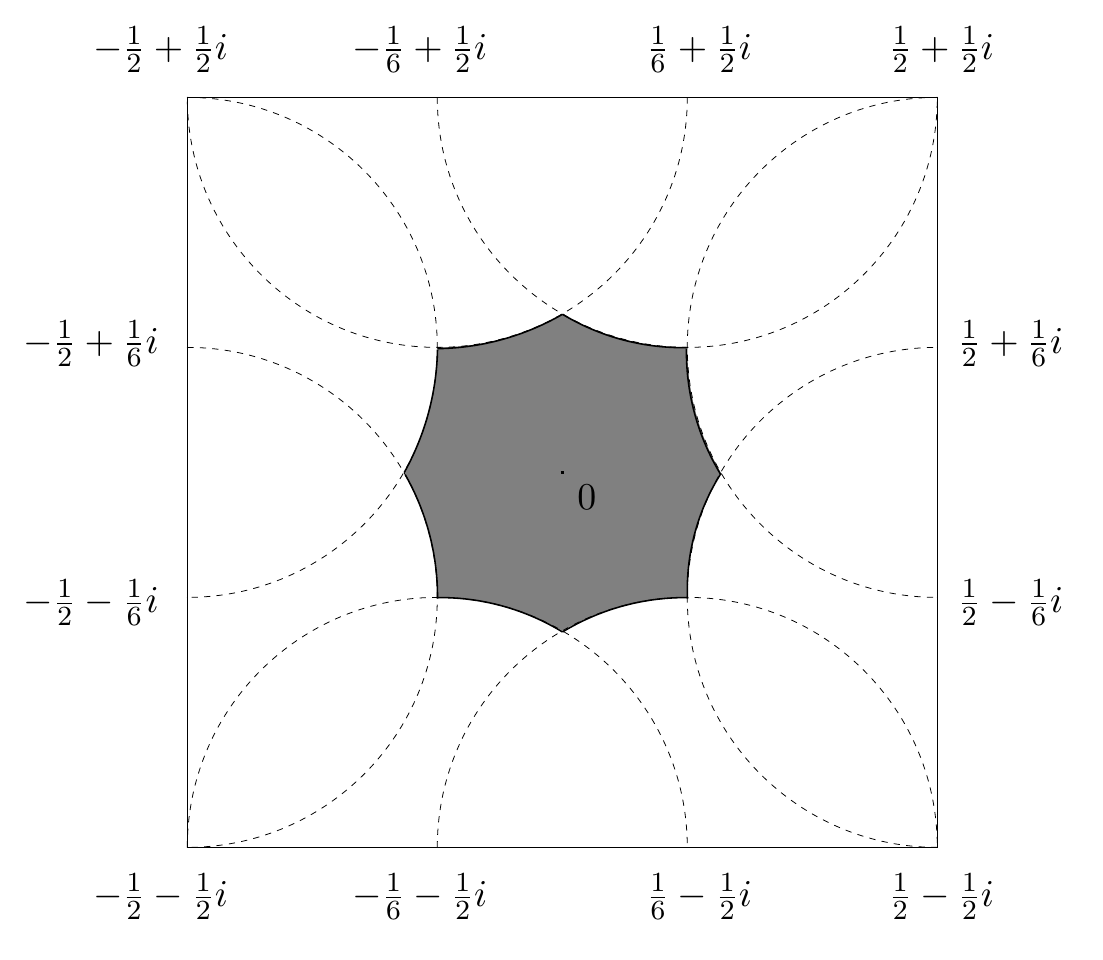}

\caption{The set $R$.}

\end{figure}

We now define a 
choice function for a continued fraction algorithm as follows: let
$z\in \C^*$ and $ 
z_0\in \frak G$ be such that $z-z_0\in S$. If $z-z_0\notin R$ we define
$f(z)$ to be the Gaussian integer nearest to $z$, and if $z-z_0 \in R$
we choose  
$f(z)$ to be the nearest odd Gaussian integer; the choice for $z_0$ as
above may not be  
unique, but the condition $z-z_0\notin R$ does not depend on the specific choice; the 
choices of the nearest integer or odd integer are not unique for certain $z$ but we
shall assume that these are fixed with some convention - this will not
play any role  
in the sequel. We shall call the algorithm corresponding to the function $f$ as above
the {\it PPOI algorithm}, using as name the acronym for ``partially preferring odd integers''.

Let $C$ denote the set of four corner points of $S$, namely $\frac 12 (\pm 1 \pm i)$, and 
for each $c\in C$ let $T(c)$ denote the subset of $S$ bounded by the triangle formed by
the other three vertices; the boundary triangle is included in $T(c)$. Also for each $c\in C$
let $\Phi (c)= c+(R\cap T(c))= \{c+z \mid z\in R\cap T(c)\}$. It can be seen that the 
fundamental set of the PPOI algorithm is given by 
$$\Phi = S \cup \bigcup_{c\in C} \Phi (c).$$ 
Furthermore, if $a\in \frak G$ is an odd integer then $f^{-1} (\{a\})$ contains all interior 
points of $a+\Phi$ and if $a$ is even $f^{-1} (\{a\})$ is disjoint from $\cup_{c\in C} \Phi (c)$
(and contains the interior of its complement in $\Phi$).

\begin{figure}[h]
\centering

\includegraphics[scale=0.8]{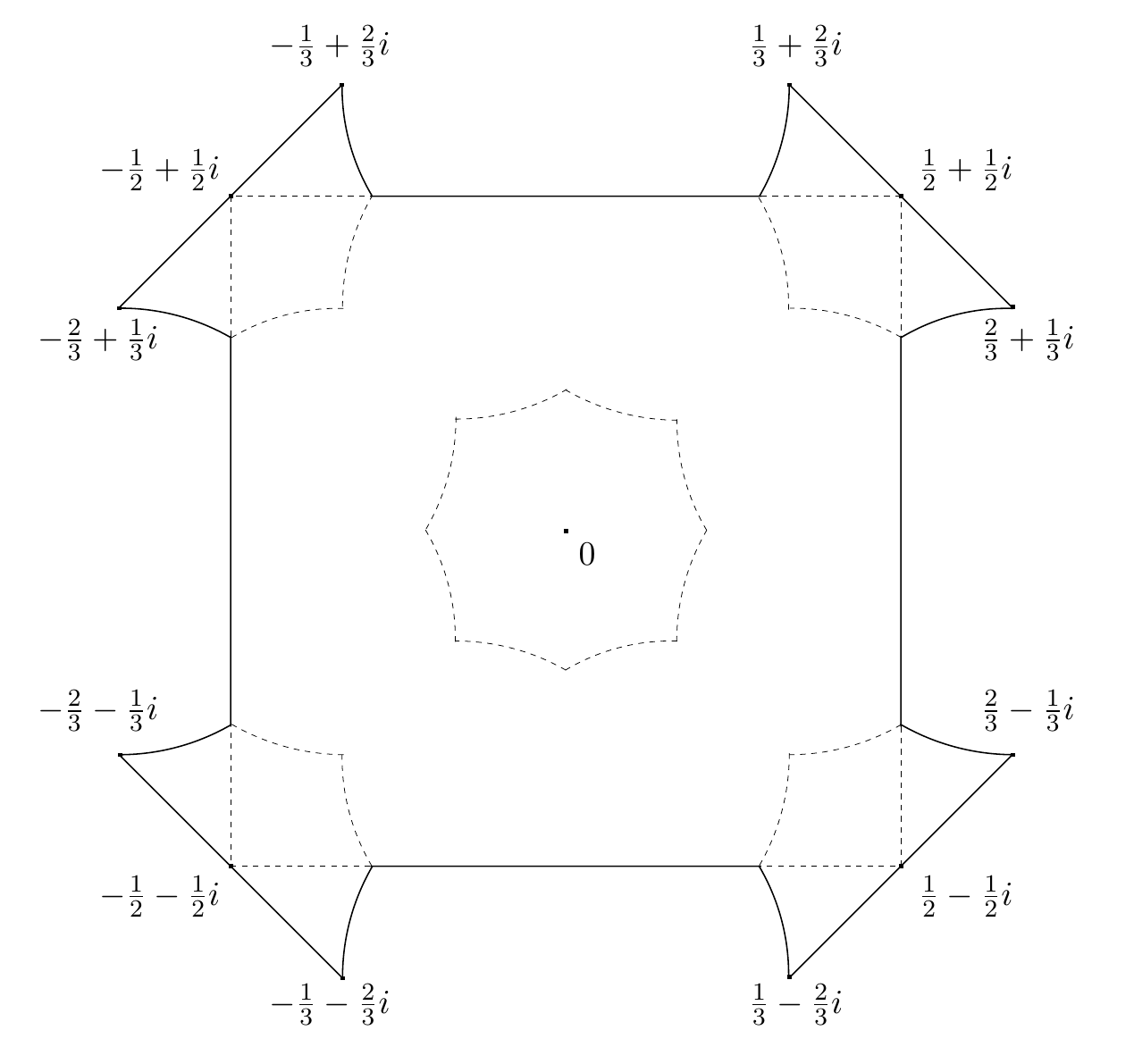}

\caption{The fundamental set $\Phi$.}
\end{figure}

\begin{proposition}\label{PPOI}
Let $z\in \C'$ and $\{a_n\}$ be the associated sequence of partial quotients with 
respect to the PPOI algorithm. Then $\{a_n\}$ satisfies Condition (H$'$) and Condition (A$'$). 

\end{proposition}

\proof We have to show that for $n\geq 1$ if $|a_{n+1}|=\sqrt 2$ then $|\bar a_n +a_{n+1}|\geq 
2$ and if $|a_{n+1}|\leq  2$ and $|a_{n}|=\sqrt 2$ then $|a_na_{n+1} +1|\geq |a_{n+1}|+1$. By 
symmetry considerations it suffices to show that if $a_{n+1}=1+i$ then $a_n$  is not one of
$-2, -1+i, 2i, -2+i, -1+2i, -2+2i$ and also not one of $-1+i, -1-i$ (correspondingly to the two 
conditions respectively), if $a_{n+1}=2$ then $a_n$  is not one of
$-2, -1+i,  -1-i$  and if $a_{n+1}=2i$ then $a_n$  is not one of $2i, 1+i,   -1+i$.  

Let $\{z_n\}$ denote that iteration sequence corresponding to $z$, with respect to the 
PPOI algorithm and let $n\geq 1$ be such that  $a_{n+1}=1+i$. Then $z_{n+1}\in (1+i+S)\cap 
\Phi^{-1}$. It can be seen that the latter is the part of $1+i+S$ in the complement of
$B(2), B(2i)$ and  $B(\frac 12 (1+i), \frac 1{\sqrt 2})$. We deduce from this  that 
$z_{n+1}^{-1}$, which is the same as $z_n-a_n$, is contained in $\Phi (\frac 12 (1-i))$. 
Since $f^{-1}(a)$ is disjoint from $\Phi (c)$ for all $c\in C$ when $a$ is even, this 
implies that $a_n$ is odd. Thus $a_n$ does not take any of the values listed above,
except possibly $-2+i$ or $-1+2i$. But  we see directly that if $a$ is one  of these
then $(a+S) \cap \Phi^{-1}$ is disjoint from $\Phi (\frac 12 (1-i))$. Thus $a_n$
does not take any of the values as above. 

Now suppose that $a_{n+1}=2$. We have $z_{n+1}\in (2+S)\cap 
\Phi^{-1}$, which may be seen to be contained in $(2+S)\backslash B(1)$. It can be readily 
deduced 
from this that $z_{n+1}^{-1}$, or equivalently $z_n-a_n$, is contained in $B(\frac 12 (1-i), 
\frac 1{\sqrt 2})$ as also $B(\frac 12 (1+i), \frac 1{\sqrt 2})$. On the other hand $(-1+i+S)\cap 
\Phi^{-1}$ and $(-1-i+S)\cap \Phi^{-1}$ are seen to be disjoint from $B(\frac 12 (1-i), 
\frac 1{\sqrt 2})$ and $B(\frac 12 (1+i), \frac 1{\sqrt 2})$ respectively. Thus $a_n$ is not 
one of $-1+i$ or $-1-i$. Similarly, observing that $z_{n+1}^{-1}$ is contained in $B(1)$ 
while  $(-2+S)\cap \Phi^{-1}$ is disjoint from $-2+B(1)=B(-1)$ we conclude that $a_n$
is not $-2$.  

A similar argument also shows that if $a_{n+1}=2i$ then $a_n$ is not one of 
$2i, 1+i,   -1+i$. This completes the proof. \hfill $\Box$

\begin{remark}
{\rm In general if $\{a_n\}_{n=0}^\infty$ is a sequence defining a continued 
fraction expansion of
a $z\in \C'$, we may not have $|z-a_0|\leq 1$. For example, let
$z'\in B(0,1)\cap B(1+i,1)\cap \C'$ and let $\{b_n\}$ be a continued fraction
expansion of $z'$  with $b_0=0$; if $z'$ is such that $0$ is the
nearest Guassian integer then $\{b_n\}$ may be taken to be the
sequence with respect to the Hurwitz algorithm - otherwise we may
use an iterated sequence construciton for this. Then for
$z=(1+i-z')^{-1}$ we see that 
while $|z|>1$, the sequence $\{a_n\}$ with $a_0=0$, $a_1=1+i$ and 
$a_n=b_{n-1}$ for all $n\geq 2$ defines a continued fraction expansion
of $z$, and we have $|z-a_0|>1$. It may also be noted that when  
$\{a_n\}_{n=0}^\infty$ defines a continued fraction expansion of
a $z\in \C'$, while for all $m\in \N$ the sequence 
$\{a_n\}_{n=m}^\infty$ defines a continued fraction expansion of a
$z_m\in \C'$, the sequence $\{z_m\}$ may not be an iteration sequence
of $z$, since the $z_m$ may not be in $\bar B(0,1)$. We shall see
below, for sequences $\{a_n\}$ satisfying 
Condition~(H$'$) $|\omega (\{a_n\})-a_0|\leq 1$ and $\omega
(\{a_n\}_{n=m}^\infty)$ is an iteration sequence of $\omega
(\{a_n\})$. 
}
\end{remark}
 
As before let $\sigma_y$ denote the map defined by $\sigma_y(z)=-\bar z$
for all $z\in \C$.  

\begin{lemma}\label{finite}
Let $\{a_n\}_{n=0}^m$ be  a finite sequence with $a_n\in \frak G$ for all $n\geq 0$ and 
$|a_n|>1$ for $n=1,\dots , m$. Let $\{p_n\}_{n=0}^m$ and $\{q_n\}_{n=0}^m$ be 
the corresponding $\cal Q$-pair. Suppose also that  for all $n\geq 1$
with $|a_n|=\sqrt 2$ the following holds: if for some $k$, with $1\leq
k\leq m-n$,  we have  $a_{n+j}= 2\sigma_y^j( a_n)$ for all $j\in \{1,\dots, 
k-1\}$ (empty set if $k=1$) and $|a_{n+k}-\sigma_y^{k} (a_n)|< 2$ then $a_{n+k}=
2\sigma_y^k( a_n)$. Then $\displaystyle{|\frac {p_m}{q_m}-a_0|< 1}$.
\end{lemma}

\proof  Let $\zeta=p_m/q_m$ and let $\zeta_0,
\zeta_1, \dots , \zeta_m$ be defined by $\zeta_0=\zeta$, and
$\zeta_{n}=(\zeta_{n-1}-a_{n-1})^{-1}$ for $n=1,\dots, m$. Then we see
that $\zeta_m=a_m$. Suppose that $|\zeta-a_0|\geq 1$.  Then we have 
$|\zeta_m|=|a_m|>1$ while on the other hand 
$|\zeta_1|=|\zeta-a_0|^{-1}\leq 1$. Hence there exists $1\leq k\leq m-1$ such that
$|\zeta_k|\leq 1$ and $|\zeta_{k+n}|> 1$ for $n=1, \dots , m-k$. 
For all $n\leq m-k$ let $z_n=\zeta_{k+n}$ and $b_n=a_{n+k}$. Then we have 
$|z_0|\leq 1$ and 
$|z_n|> 1$ for $n=1,\dots , m-k$. Also $z_n-z_{n+1}^{-1}=\zeta_{n+k}-\zeta_{n+k+1}^{-1}
=a_{n+k}=b_n$. Applying Lemma~\ref{simple} and using the condition as
in the hypothesis  we get that $|\zeta_m|=|z_{m-k}|\leq \sqrt 2+1$
and $|a_m|=|b_{m-k}|=2\sqrt 2$. But this is a contradiction since $\zeta_m=a_m$. 
Therefore $|\displaystyle{\frac {p_m}{q_m}-a_0|< 1}$. \hfill $\Box$

We can now conclude the following. 

\begin{theorem}\label{represent}
Every irrational complex number  $z$ admits a continued fraction expansion 
$\{a_n\}_{n=0}^\infty$, $a_n\in \frak G$ for all $n$, satisfying the
following conditions:  

i) $|a_n|>1$ for all $n\geq 1$,

ii) $|a_{n}+\bar a_{n+1}|\geq 2$ for all $n\geq 1$ such that
$|a_{n+1}|=\sqrt 2$, and  

iii) $|a_na_{n+1}+1|\geq |a_{n+1}|+1$ for all $n\geq 1$ such that
$|a_{n}|=\sqrt 2$  and $|a_{n+1}|\leq  2$. 

Conversely, if $\{a_n\}_{n=0}^\infty$, $a_n\in \frak G$, is a sequence 
satisfying 
conditions~(i), (ii) and (iii), and 
$\{p_n\}$, $\{q_n\}$ is the corresponding $\cal Q$-pair, then the
following holds:

a) $q_n\neq 0$ for all $n\geq 0$ and $p_n/q_n$ converges to an 
irrational complex number $z$ as $n\to \infty$, (so $\{a_n\}$ defines 
a continued fraction expansion of $z$);

b) $|z-a_0|\leq 1$ and $\{z_n\}$ is an iteration sequence for $z$;

c) $z$ is a quadratic surd if and only if $\{a_n\mid n\in \N\}$ is 
finite. 
\end{theorem}

\proof The first assertion follows from Proposition~\ref{PPOI}. 
In  the converse part, assertion~(a) follows from Proposition~\ref{growth2}
and Theorem~\ref{cgce}. Lemma~\ref{finite} implies by 
condition~(ii)  that $|z-a_0|\leq 1$, and applying 
the argument successively to $\omega (\{a_n\}_{n=k}^\infty\}$,
$k=1,2, \dots$ we see that  $\{z_n\}$ is an iteration sequence for
$z$, and thus (b) is proved. Assertion~(c) follows from
Proposition~\ref{dir} and Corollary~\ref{lagr-iter}. \hfill $\Box$

\section{Associated fractional transformations}

As with the usual continued fractions we can associate with the complex continued 
fractions  transformations of the projective space, the complex projective
space in this instance.

We denote by $\C^2$ the two-dimensional complex vector space, viewed as
the space of 2-rowed columns with complex number entries. 
We denote by $e_1$ and $e_2$ the  vectors $ \left (   
\begin{array}{c}  
      1\\
     0\\ 
      \end{array} \right) $ and 
      $ \left (   
\begin{array}{c}  
      0\\
     1\\ 
      \end{array} \right),$
      in $\C^2$, respectively. For all $z\in \C$ we denote by $v_z$ the vector $ze_1+e_2$. 

We denote by 
$\P$ the corresponding one-dimensional complex projective space, consisting 
of complex lines (one-dimensional vector subspaces)  in $\C^2$ and by 
$\pi :  \C \to \P$ the map defined by  setting, for all $z\in \C$,
$\pi (z)$ to be the complex line containing $ v_z$.

We denote by $SL(2,\C)$ the group of $2\times 2 $ matrices with complex entries,
with determinant 1. The natural linear action of $SL(2,\C)$ on $\C^2$ induces an 
action of the group on $\P$. We shall be using this action freely in the sequel, with 
no further mention. We denote by $SL(2,\frak G)$ the subgroup of $SL(2,\C)$ consisting 
of all its elements with entries in $\frak G$. 

When $\C$ is realised as a subset of $\P$ via the map $\pi$ as above, the 
subset $\C'$ (as before, consisting of complex numbers which are
not rational) is readily seen to be invariant under the action of $SL(2,\frak G)$; the 
action of the element $\left (   
      \begin{array}{cc}  
      p &  q\\
      r& s\\ 
      \end{array} \right) \in  SL(2,\frak G)$ on $\C'$ is then given by 
      $z\mapsto \displaystyle{\frac{pz+q}{rz+s}}$ for all $z\in \C'$. 

We denote by $\rho $ the matrix defined by 
$$\rho =   \left (   
      \begin{array}{cc}  
      0 &  -1\\
      1& 0\\ 
      \end{array} \right).
      $$
For any $g\in SL(2,\C)$ we denote by $^t\!g$ the transpose of the matrix $g$. It 
can be verified that for any $g\in SL(2,\C)$, $^t\!g\rho g =\rho $; this identity will be used 
crucially below. 
For any $a \in \frak G$ we  denote by $ U_a$ and $L_a$ the matrices  defined by 
$$U_a=   \left (   
      \begin{array}{cc}  
      1 &  a\\
      0& 1\\ 
      \end{array} \right) \hbox{\rm and  }
      L_a=   \left (   
      \begin{array}{cc}  
      1 &  0\\
      a & 1\\ 
      \end{array} \right).
      $$
     
By a {\it block} $\beta $ of Gaussian integers we mean a finite $l$-tuple 
$\beta= (b_1, \dots , b_l)$,  where $l\geq 1$ and 
$b_1, \dots , b_l \in \frak G$; the $l$ is called the {\it length} of the block and 
will be denoted by $l(\beta)$.
For a block $\beta = (b_1, \dots , b_l)$ we define 
$$g(\beta) = U_{b_1}L_{b_2}\cdots U_{b_{l-1}}L_{b_{l}} \ \hbox{ \rm or  } \ 
U_{b_1}L_{b_2}\cdots L_{b_{l-1}}U_{b_{l}},$$
according to whether $l$ is even or odd, respectively.    

A block $\beta = (b_1, \dots , b_l)$ is said to {\it occur} in a
sequence $\{a_n\}$ in 
$\frak G$, starting  at $k\geq 0$ if  $a_k=b_1$, $a_{k+1}=b_2, ...,
a_{k+l-1}=b_l$, and it is said to be an {\it initial block} of $\{a_n\}$ if it occurs starting at $0$. 

\begin{remark}\label{cf}
{\rm Let  $\{a_n\}$  be a sequence from $\Omega_E^{(\theta)}$, $\theta >1$,
and $\beta = (b_1, \dots , b_l)$ be an initial block of  $\{a_n\}$. Let  
$\{p_n\}$, $\{q_n\}$ be the $\cal Q$-pair corresponding to $\{a_n\}$.  Then it is straightforward 
to verify, say by induction, that  for $r=1$ and $2$, 
$g(\beta)(e_r)=  p_me_1+q_me_2$, where  $m$ is the largest odd integer less than 
$l(\beta)$ if $r=1$ and the largest even integer  
less than $l(\beta)$ if $r=2$. 
}
\end{remark}

For $\theta >1$, a sequence $\{\beta_i\}$ of blocks of Gaussian integers is said to
{\it $\theta$-represent} a number  $z\in \C'$ if $l(\beta_i)\to \infty$ and for 
each $i$ there exists 
a sequence $\{a_n\}$ (depending on $i$) in $\Omega_E^{(\theta)}$ such that $\omega (\{a_n\})=z$ and 
$\beta_i$ is an initial block of $\{a_n\}$.

\begin{theorem} \label{cfl} 
Let $\zeta \in \C'$ and $\{\beta_i\}$ be a sequence of blocks of Gaussian integers 
$\theta$-representing $\zeta$.  Then as $i\to \infty$,
$$g(\beta_i)(\pi (z)) \to \pi (\zeta) \hbox{ \rm for all } z \hbox{ \rm with }  |z|>1,$$
uniformly over $\{z \in \C \mid   |z|\geq 1+\epsilon \}$, for all $\epsilon >0$.  
\end{theorem} 

\proof It would suffice to prove the assertion 
in the Theorem~separately with all $l(\beta_i)$ being odd and 
all being even. We shall give the argument when all $l(\beta_i)$ are even; a similar 
argument works in the other case. 
Let $z\in \C$ with $|z|>1$ be given. For each $i$ let $j_i$ and $k_i$ be the 
largest odd and even integers less than $l(\beta_i)$, respectively;
as  $l(\beta_i)$ are even we have $j_i=l(\beta_i)-1$ and $k_i=j_i-1$
for all $i$. Then by Remark \ref{cf} we have  
$$g(\beta_i)(\pi (z))=\pi (\frac {zp_{j_i}+p_{k_i}}{zq_{j_i}+q_{k_i}})= 
\pi (\frac{p_{j_i}}{q_{j_i}}\frac {z+(p_{k_i}/p_{j_i})}{z+(q_{k_i}/q_{j_i})}).$$
Since $\{\beta_i\}$ is a sequence of blocks  $\theta$-representing $\zeta$, 
by Corollary \ref{commoncf} we get that $\{\frac{p_{j_i}}{q_{j_i}}\}$ and 
$\{\frac{p_{k_i}}{q_{k_i}}\}$ converge to $\zeta$. Since $\{a_n\}$ belongs to 
$\Omega_E^{(\theta)}$ and $k_i<j_i$ we have 
 $|q_{k_i}|<|q_{j_i}| $ for all $i$. As $|z|>1$ this shows that  $ \{z+(q_{k_i}/q_{j_i})\}$ 
 is a bounded sequence,
 bounded away from $0$; moreover,  whenever  a subsequence of 
 $ \{z+(q_{{k_i}}/q_{j_i})\}$ converges (to a nonzero 
 number), the corresponding subsequence of $ \{z+(p_{{k_i}}/p_{j_i})\}$ also converges 
 to the same number, since $\frac {p_{{k_i}}/p_{j_i}}{q_{{k_i}}/q_{j_i}}=\frac{p_{k_i}}{q_{k_i}}
 \frac{q_{j_i}}{p_{j_i}}\to \zeta \cdot \zeta^{-1}=1$.  Thus 
 we get that $g(\beta_i)(\pi (z)) \to \pi (\zeta)$ as $i\to
 \infty$. Moreover,  
 clearly, for any $\epsilon >0$  the convergence is uniform over $\{z \in \C \mid   |z| \geq 1
 +\epsilon \}$.  \hfill $\Box$ 

\medskip
We next note also the following simple fact which can be proved along the lines of the 
proof of Lemma~2.4 of \cite{DN}; we omit the details. 

\begin{proposition}  \label{standard} 
Let $\{g_i\}$ be an unbounded sequence in $SL (2,\C)$. Suppose that 
$\varphi, \psi$ in $\P$ such that the sets $\{\xi \in \P\mid g_i(\xi) \to \varphi \}$ 
and $\{\xi \in \P\mid \, ^t\!g_i(\xi) \to \psi \}$ have at least two elements each, and 
that $g_i(\rho (\psi ))$ converges. Then $g_i(\xi) \to \varphi $ for all $\xi \neq \rho  (\psi)$.
\end{proposition}

Proposition~\ref{finite} readily implies the following.

\begin{corollary}\label{transpose}
Let  $\beta =(b_1, \dots, b_l)$ be a block of Gaussian integers, with $l(\beta)$ 
even, occurring in  a sequence 
$\{a_n\}$ satisfying Condition (H), starting at   $k\geq 0$.  If $z_r\in \C $, $r=1,2$,
be the numbers defined by  $\pi (z_r)= \, ^t\!g(\beta)(\pi (e_r))$. Then ${|z_r-b_{l}|<1}$. 
\end{corollary}

We now describe the main technical result towards construction of ``generic points'',
after introducing some more notation and terminology. For a block $\beta =
(b_1, \dots, b_l)$ we denote
by $\tilde \beta$ the block $(b_l, \dots, b_1)$, with the entries arranged in the reverse
order. Given two blocks $\beta =(b_1, \dots, b_l)$ and $\beta' =(b'_1, \dots, b'_{l'})$,
we denote by $\beta \beta'$ the block $(b_1, \dots, b_l, b'_1, \dots, b'_{l'})$. A block
$\beta' =(b'_1, \dots, b'_{l'})$ will be called the {\it end-block} of $\beta =(b_1, \dots, b_l)$
if $b'_1=b_{l-l'+1},  \dots, b'_{l'}=b_{l}$. For a block $\beta =(b_1, \dots, b_l)$, $b_1$ and
$b_l$ will be called the first and last entries respectively. 

\begin{proposition}\label{prop:limits}
Let  $\{x_n\}$ be a sequence in $ \Omega_E^{(\theta)}$, $\theta >1$, and let $\xi =
\omega (\{x_n\})$. Let $\alpha, \beta \in \C'$ and $\{\alpha_i\}$ and  $\{\beta_i\}$ be 
sequences of blocks of Gaussian integers $\theta$-representing 
$\alpha$ and $\beta$ respectively. Suppose that the following conditions are satisfied: 

i) for all $i$, $\alpha_i$ and $\beta_i$ are of even length and their first entries have
 absolute value greater than $1$, and the last entry of $\alpha_i$ has absolute value 
greater than $2$;

ii) there exists a sequence $\{s_i\}$ of even integers such that for all $i$ the block 
$\tilde \alpha_i \beta_i$ occurs in $\{x_n\}$ starting at $s_i$, $s_{i+1}> s_i +
l(\alpha_i)+l(\beta_i)$,  $|x_{s_i -1}|>2$ and $|x_{s_i +l(\alpha_i)+l(\beta_i)}|>2$.

\smallskip
\noindent For all $i$ let $\xi_i=(x_0, \dots , x_{s_i -1})$ and $g_i=g(\alpha_i)^t\!g (\xi_i)$. 
Then, as $i\to \infty$,
$$g_i (\pi (z)) \to \pi (\alpha) \hbox{ \rm for all } z\neq \rho (\pi (\xi)), \hbox{ \rm and} \ 
g_i (\rho (\pi (\xi))) \to \rho (\pi (\beta)).$$ 
\end{proposition}

\proof Let $i \in \N$ and  for $r=1,2$ let $z_r$ be defined by 
$\pi (z_r)= {^t\!g} (\xi_i)(\pi (e_r))$. Since $|x_{s_{i -1}}|>2$, and 
hence at least $\sqrt 5$, by Corollary \ref{transpose} it follows that $|z_r|\geq \sqrt 5 -1>1$.
Since this holds for all $i$ by Proposition~\ref{cfl} we get that, as $i\to \infty$, 
$g_i (\pi (e_r))= g(\alpha_i)(\pi (z_r))\to \pi (\alpha)$, for $r=1$ and $2$. Similarly, 
under the conditions in the hypothesis we get  $^t\!g_i (\pi (e_r)) \to \pi (\xi) $ for $r=1,2$.

We next show that $g_i (\rho (\pi (\xi)))\to \rho (\pi (\beta))
$. Then by Proposition~\ref{standard} this would imply the conclusion
as stated in the present proposition.  
Let $1\leq i < j$. Then from the conditions in the hypothesis we see that $\xi_j$ may be
written as  $\xi_j= \xi_i\tilde \alpha_i \beta_i \varphi_{i,j}$, where $\varphi_{i,j}$ is the 
end-block of $\xi_j$  of length $s_{i+1}-s_i-l(\alpha_i)-l(\beta_i)>0$.
Therefore $g (\xi_j)=g(\xi_i) ^t\!g(\alpha_i)g(\beta_i)g(\varphi_{i,j})$,
and hence $$g_i \rho  g (\xi_j)= g(\alpha_i) ^t\!g (\xi_i)\rho g(\xi_i) 
^t\!g(\alpha_i)g(\beta_i)g(\varphi_{i,j})=\rho g(\beta_i)g(\varphi_{i,j}),$$
using (twice) the fact  that $^t\!g\rho g=\rho $ for all $g\in SL(2,\C)$. 
We note that $s_i +l(\alpha_i)
+l(\beta_i)$ are even and $x_{s_i +l(\alpha_i)+l(\beta_i)}$ is the first entry of 
$\varphi_{i,j}$. Since  $|x_{s_i +l(\alpha_i)+l(\beta_i)}|\geq \sqrt 5$, arguing as in 
the first part we conclude that $\rho g(\beta_i)g(\varphi_{i,j}) $ 
converges to 
$\rho (\pi (\beta))$, as $i$ and $j$ tend to infinity, for $r=1,2$. Hence $g_i \rho  g (\xi_j)(\pi(e_r))$,
$r=1,2$, converge to  $\rho (\pi (\beta))$, as $i$ and $j$ tend to infinity. Since any limit 
point of the sequence $\{g_i (\rho (\pi (\xi)))\}$ is a limit point of $g_i \rho  g (\xi_j)(e_r)$ as 
$i$ and $j$ tend to $\infty$, this shows that 
$g_i (\rho (\pi (\xi)))\to \rho (\pi (\beta)) $. This completes the proof of the Proposition.
\hfill $\Box$

\section{Application to values of quadratic forms}

In this section we apply the results on complex continued fractions to the study of 
values of quadratic forms 
with complex coefficients over the set of pairs of Gaussian integers. It was
proved by G. A. Margulis in response to a long standing conjecture due to 
A. Oppenheim that for any nondegenerate indefinite real quadratic form $Q$
in $n\geq 3$ variables which is not a scalar multiple of a rational quadratic
form, the set of values of $Q$ on the set of $n$-tuples of integers forms a dense 
subset of $\R$. This was generalised by A. Borel and G. Prasad to 
quadratic forms over other locally compact fields \cite{BP}, and in particular
it is known that for a nondegenerate complex quadratic form in $n\geq 3$ variables
which is not 
a scalar multiple of a rational form, the values on the set of $n$-tuples of Gaussian
integers are dense in $\C$. The reader is referred to \cite{D-encl} and 
\cite{M-survey} for a discussion on the overall theme, including Ratner's 
pioneering work on dynamics of unipotent flows on homogeneous spaces
to which the results are related. The case of binary quadratic forms, viz. 
$n=2$, does not fall in the ambit of the theme, and in fact the corresponding
statement is not true for all irrational indefinite binary quadratic forms. For example
if $\lambda $ is a badly approximable number (here we may take a real number
whose usual partial quotients are bounded) then for the values of the quadratic 
form $Q(z_1, z_2)=(z_1-\lambda z_2)z_2$, it  may be seen that $0$ is an isolated 
point of the set of values of $Q$ over the set of pairs of Gaussian integers. 

In \cite{DN} we described a  class of real binary quadratic 
forms for which the values on the set of pairs of integers is dense in $\R$;
this is related to a classical result of E. Artin on orbits of the
geodesic flow associated  
to the modular surface; in \cite{DN}  Artin's result was 
strengthened, and an application was made to  values of the binary
quadratic forms over  
the set of pairs of integers, and also over the set of pairs of positive integers. 

In this section, using  complex 
continued fractions we describe an analogous class of complex binary
quadratic forms for which the set of values on pairs of Gaussian integers
is a dense subset of $\C$. 

 For 
$ u=\left (   
\begin{array}{c}  
      \alpha_1\\
     \alpha_2\\ 
      \end{array} \right) $ and 
      $ v=\left (   
\begin{array}{c}  
      \beta_1\\
     \beta_2\\ 
      \end{array} \right), $
we denote by $Q_{u,v}$ the quadratic form on $\C^2$ defined by 
$$Q_{u,v}(z_1e_1+z_2e_2)=(\alpha_1z_1+\alpha_2z_2)(\beta_1z_1+\beta_2z_2) 
\hbox{ \rm for all } z_1,z_2 \in \C.$$

We denote by $\eta : \C^2\backslash \{0\} \to \P$ the canonical quotient map. We note 
the following fact which can be proved  in the same way as Proposition~5.1 of \cite{DN}. 

\begin{proposition}\label{limits}
Let $u,v, u',v' \in \C^2\backslash \{0\}$ and suppose that there exists a sequence $\{g_i\}$ in 
$SL(2,\C)$ such that 
$g_i(\eta (u))\to \eta (u')$ and $g_i(\eta (v))\to \eta (v')$, as $i\to \infty$, in $\P$. Let $E$ be 
a subset of $\C^2$ such that
$^t\!g_i (E)\subseteq E$ for all $i$. Then the closure of the set $Q_{u,v}(E)$ in $\C$  contains 
$Q_{u',v'}(E)$.
\end{proposition}

We denote by ${\frak G}^2$ the subgroup of $\C^2$ consisting of pairs of Gaussian integers,
namely $\{ae_1+be_2\mid a,b \in \frak G\}$. We note that ${\frak G}^2$ is invariant under the 
linear action of $SL(2, \frak G)$ on $\C^2$. 

\begin{remark}\label{orbitconnection}
{\rm Consider the quadratic form $Q=Q_{u,v}$ where $u$ and $v$ are linearly independent. 
We write $u$ and $v$ as $g^{-1}(e_1),g^{-1}(e_2)$  where $g\in GL(2,\C)$, so $Q=Q_{g^{-1}
(e_1),g^{-1}(e_2)}$.  Replacing $Q$
by a scalar multiple we may assume (in studying the question of values over ${\frak G}^2$
being dense) $g$ to be in $SL(2,\C)$.  Let $D$ be the subgroup of $SL(2,\C)$ consisting of
diagonal matrices. Then we see that for any $d\in D$, $\eta (u)$ and
$\eta (v)$ are fixed  
under the action of $g^{-1}dg$. If there exist sequences $\{\gamma_i\}$ in $SL(2,{\frak G})$ 
and $\{d_i\}$ in $D$ such that $\{\gamma_ig^{-1}d_ig\}$ is dense in $SL(2,\C)$ then by 
Proposition~\ref{limits} it follows that $Q_{u,v}({\frak G}^2)$ is dense in $\C$. The condition 
involved is equivalent to that the orbit of the coset of $g$ in $SL(2,\C)/SL(2,{\frak G})$ under the 
action of $D$ is dense in $SL(2,\C)/SL(2,{\frak G})$. It is well known that $SL(2,\frak G)$ is 
a lattice in $SL(2,\C)$, viz. $SL(2,\C)/SL(2,{\frak G})$
admits a finite measure invariant under the action of $SL(2,\C)$ (see \cite{Bor}, \S 13 for a 
general result in this respect). This implies in particular  that the action of $D$ on
$SL(2,\C)/SL(2,{\frak G})$ is ergodic (with respect to the $SL(2,\C)$-invariant measure); see
\cite{BM}, Chapter III, for instance. Hence for {\it almost all} $g$ the above condition holds and 
consequently the values of $Q=Q_{g^{-1}(e_1),g^{-1}(e_2)}$ over ${\frak G}^2$ are 
dense in $\C$. However, this does not enable describing any specific quadratic form 
for which the conclusion holds. The point about the application here is that we obtain a set of 
forms which may be explicitly described, in terms of continued fraction expansions, for
which this property holds.

}
\end{remark}

We next define the notion of a generic complex number. In this respect
we shall restrict  
to sequences arising as continued fraction expansion with respect to
the Hurwitz algorithm;   
the procedure can be extended to other algorithms satisfying suitable
conditions, but in 
the present instance we shall content ourselves considering only the
Hurwitz algorithm.  

We shall assume that the algorithm is determined completely, fixing 
a convention  
for tie-breaking when there are more than one Gaussian integers
at the minimum distance, so that  
we can talk unambiguously of the continued fraction expansion of
any number; for brevity we shall refer to the corresponding sequence
of partial quotients as the {\it Hurwitz
expansion} of the number.  
Let $\Omega_H$ denote the set of sequences $\{a_n\}$ such that
$\{a_n\}$  is the Hurwitz  expansion of a $z\in \C'$. As before let
$S=\{z=x+iy\mid |x|< \frac 12,  |y|< \frac 12\}$. 
For a block $\beta=(b_1, \dots, b_l)$, with $b_i\in  \frak
G$ and $|b_i|>1$ for $i=1, \dots ,l$,  we denote by $S_\beta$ the subset
of $S$ consisting of  
all $z\in S$  with Hurwitz expansion $\{a_n\}$ such that
$(a_1, \dots a_l)=\beta$   (note that since $z\in S$, $a_0=0$). 
We say that a block $\beta=(b_1, \dots, b_l)$ is {\it admissible} if 
$S_\beta$ (is nonempty and) has nonempty interior. We note that
$S_\beta$ is a set 
with a boundary consisting of finitely many arc segments (circular or 
straight line) and the
interior, if nonempty, is dense in $S_\beta$. (It seems plausible that
$S_\beta$ has nonempty interior whenever it is nonempty, but it does
not seem to be easy to ascertain this - as this will not play a role
in the sequel we shall leave it at that.)

\smallskip
We say that $\{a_n\}\in \Omega_H$ is {\it generic} if every
admissible block occurs in  
$\{a_n\}$, and $z\in \C'$ is said to be {\it generic} if its Hurwitz 
expansion is a generic sequence. 

\begin{remark}\label{follow}
{\rm 
We note the following facts which may be verified from inspection of the iteration 
sequences corresponding to the Hurwitz expansions.   
If $\beta=(b_1,\dots , b_l)$ and $\beta'=(b_1',\dots , b_{l'}')$ are admissible blocks and 
$|b_l |\geq 2\sqrt 2$ then $(b_1,\dots, b_l, b_1',\dots , b_l')$
is an admissible block. If $|b_l |<2\sqrt 2$ then in general
$(b_1,\dots , b_l, b)$ need not be 
admissible for $b\in \frak G$, $|b|>1$. However, given any
admissible block  
$\beta=(b_1,\dots , b_l)$ and $b'\in \frak G$ with $|b'| \geq 2\sqrt
2$ there exists $b\in \frak G$ such that $|b|=|b'|$  
and $(b_1,\dots, b_l, b)$ is admissible. This shows in particular that given admissible blocks
$\beta$ and $\beta'$ there exists an admissible block
$(x_1,\dots, x_m)$ such that $(x_1,\dots, x_l)=\beta $ and $(x_{m-l'+1},\dots, x_m)=\beta'$,
where $l$ and $l'$ are the lengths of $\beta $ and $\beta'$
respectively, and furthermore $m$ may be chosen to be at most $l+l'+1$. 
Since the set of admissible blocks is countable, this also shows that we can construct
generic complex numbers, by building up their Hurwitz expansion  in steps, so
that  each admissible block occurs at some stage.}
\end{remark}

\begin{remark}\label{midpoint}
{\rm For any block $\beta=(b_1,\dots , b_l)$ let $r_\beta =\omega
  (\{b_n\}_{n=0}^l)$ (notation as introduced \S~2), with $b_0=0$. Then
  each $r_\beta $ is a rational number from $S$; if all $b_i,
  i=1,\dots ,l$ are of absolute value greater than $2$ then
$r_\beta$ is an interior point of $S_\beta$ -otherwise it is a
boundary point of $\beta$. Let $Q$ be the set
consisting of all $r_\beta$ with $\beta$ any admissible block. Then
$Q$ is a dense subset  
of $S$  since for any nonempty open subset of $S$ we can find a block
$\theta$ such  
that $S_\theta$ is contained in it. 
}

\end{remark}

\begin{proposition}\label{Gdelta}
The set of generic complex numbers contains a dense ${\cal G}_\delta$ subset of $\C$. 
\end{proposition}

\proof  Let $\beta$ be any admissible block of length $l\geq 1$ and
let $E_\beta$ 
be the interior of the set of numbers $z\in S$ such that $\beta$
occurs in the Hurwitz 
expansion of $z$.   Consider any admissible block $\gamma
=(c_1, \dots , c_k)$. Let $\epsilon >0$ be arbitrary and let $N\in \N$
such that $N>\epsilon^{-1}$. By  
Remark~\ref{follow} there exists an admissible block $\xi=(x_1,\dots, x_m)$ 
such that $(x_1,\dots , x_k)=\gamma$, $x_{k+1}=N+2$ and $(x_{m-l+1},
\dots , x_m)=\beta$.  Then the  interior of $S_\xi$ is contained in
$E_\beta$. We note that for any $z\in S_\xi$ by Proposition~\ref{prop} we have
$|z-r_\gamma|\leq |z_{k+1}-1|^{-1}$, where $\{z_n\}$ denotes the iteration
sequence for $z$ with respect to the Hurwitz algorithm,  and we have 
$|z_{k+1}|\geq (N+2)-|z_{k+2}|^{-1}>N+1$; see  Proposition~\ref{prop}. Hence 
the distance of $r_\gamma$ from $S_\xi$ is at most $N^{-1}<\epsilon$. 
Since $\epsilon >0$ is arbitrary this shows that $r_\gamma$ is
contained in the closure of $S_\xi$. As the  interior of $S_\xi$ is 
a dense subset of $S_\xi$ contained in $E_\beta$, this shows that 
$r_\gamma$ is contained in the closure of $E_\beta$. Since this holds for
every admissible block $\gamma$ 
it means that the set $Q$ as defined in Remark~\ref{midpoint} 
is contained in $E_\beta$, and since $Q$ is dense in $S$ we get that 
$E_\beta$ is dense in $S$.  Thus for every admissible block $\beta$
there exists an open dense subset of $S$ consisting of points such
that $\beta$ occurs in their Hurwitz expansion. As there are only
countably many blocks we get that there exists a dense ${\cal
  G}_\delta$ subset of $S$ consisting of generic complex numbers. Since
translates of generic numbers by  Gaussian integers are generic it 
follows that there exists a dense ${\cal G}_\delta$ subset of $\C$
all whose elements are generic complex numbers.  \hfill $\Box$

\medskip
In the light of a result of \cite{Hen} we can conclude also the following. 

\begin{proposition}\label{fullmeasure}
With respect to the Lebesgue measure on $\C$ almost every $z$ in $\C$ is generic. 
 \end{proposition}

\proof Let $S=\{z=x+iy\mid |x|\leq \frac 12, |y|\leq \frac 12\}$. For any  
admissible block $\beta$  consider the subset, say $E_\beta$, of $S$ consisting of all $z$ 
such that $\beta$ occurs in the  continued expansion of $z$. We note that $E_\beta$ is 
invariant under the Gauss map, of $S$ into itself, associated with the Hurwitz continued fractions
(see \cite{Hen}). It is known that the Gauss map is ergodic with respect to the Lebesgue 
measure; in fact  there exists a probability measure on $S$ equivalent to the 
Lebesgue measure which is invariant with respect to the Gauss map, and
the map is mixing with respect to it; see \cite{Hen}, Theorem~5.5. Since $E_\beta$
contains an nonempty open set it is of positive measure, and hence by ergodicity
it is of measure $1$. It follows that the set of generic points in $S$, which is the intersection 
of $E_\beta$ over the countable collection of admissible blocks $\beta$, is of full measure in $S$. 
Since translates of generic points by Gaussian integers are generic,
we get that 
almost every complex number is generic.  \hfill $\Box$

\medskip

 We now note the following.

\begin{corollary}\label{dense}
Let $\xi \in \C$ be generic.
Let $\alpha, \beta \in \C$ and $\{a_n\}$ and $\{b_n\}$ their respective continued fraction 
expansions.   Suppose that $|a_n|> \sqrt 5$ 
for all $n\geq 0$,  $|b_0|\geq \sqrt 2$, and all initial blocks of $\{b_n\}$ are admissible.
Then there exists a sequence $\{g_i\}$ in  $SL(2,\frak G)$ such that 
$$g_i (\pi (z)) \to \pi (\alpha) \hbox{ \rm for all } z\neq \rho (\pi (\xi)), \hbox{ \rm and} \ 
g_i (\rho (\pi (\xi))) \to \rho (\pi (\beta)).$$
\end{corollary}

\proof We note that under the conditions on $\alpha$ and  $\beta $ in the hypothesis, 
if $\alpha_i$ and $\beta_i $ are initial blocks of length $2i$ in $\{a_n\}$ and  $\{b_n\}$
respectively, then $\tilde \alpha_i \beta_i$  are 
admissible blocks.  As $\{x_n\}$ is generic we can find a sequence $\{s_i\}$ in $\N$ 
such that the conditions as in the hypothesis 
of Proposition~\ref{prop:limits} are satisfied. The assertion as in the Corollary then
follows from the Proposition~\ref{prop:limits}. \hfill $\Box$

\begin{corollary}\label{cor:qf}
Let   $\zeta_1, \zeta_2\in \C$  and suppose that $\zeta_1$ (or $\zeta_2$) is generic. Then the values of the quadratic 
form $Q(x,y)=(x -\zeta_1y)(x -\zeta_2 y)$ over the set of Gaussian integers form a dense subset of $\C$;
namely, $\{(a - \zeta_1b)(a -\zeta_2 b)\mid a,b \in \frak G\}$ is dense in $\C$.

\end{corollary}

\proof  Let $v_r=\left (   
\begin{array}{c}  
      \zeta_r\\
     1\\ 
      \end{array} \right) $
      for $r=1,2$.  
We shall show that 
for the quadratic form $Q=Q_{\rho (v_1),\rho (v_2)}$ on $\C^2$, $Q({\frak G}^2)$ is dense in $\C$.
This is readily seen to hold  in the case when $\zeta_1=\zeta_2$, we shall 
now assume that $\zeta_1\neq \zeta_2$. We have $\eta (v_1)= \pi (\zeta_1)$,  $\eta (v_2)= \pi (\zeta_2)$ and $\eta (v_1) \neq \eta (v_2)$. 
Let $\alpha, \beta \in \C$ be such that the conditions in the statement of Corollary \ref{dense} 
are satisfied. Then by  Corollary \ref{dense} there exists a sequence $\{g_i\}$ in 
$SL(2,\frak G)$ such that
$g_i (\rho (\eta(v_1))) \to \rho (\pi (v_\beta))$ and $g_i (\eta (v_2)) \to \pi (v_\alpha)$.
 Then by Proposition~\ref{fullmeasure} the closure of 
$Q({\frak G}^2)$ in $\C$ contains $Q_{\rho (v_\beta), v_\alpha}({\frak G}^2)$, for all $\alpha, \beta \in \C$ 
as above. We note that the condition in the corollary is 
satisfied for all $\beta$ in a dense subset of  $\{z\in \C \mid |z|\geq \sqrt {5/2}\}$; it may be 
noted in this respect that the set of points for which all initial blocks of the Hurwitz continued
fraction expansion are admissible is a set of full Lebesgue measure, containing the 
complement of a countable union of analytic curves.
We fix an $\alpha$ for which the condition in  Corollary \ref{dense} is satisfied and 
consider values $Q_{\rho (v_\beta), v_\alpha}({\frak G}^2)$, for varying $\beta$ as above. It is
straightforward to see that these values contain a dense set of complex numbers. 
Hence $Q({\frak G}^2)$ is dense in $\C$. \hfill $\Box$

\begin{remark}{\rm The set of all (complex) binary quadratic forms may be viewed  
in a natural way as the $3$-dimensional vector space (over $\C$) of symmetric $2\times 2$ 
matrices, and may be considered equipped with the Lebesgue measure. From the fact 
that the generic numbers in $\C$
form a set of full Lebesgue measure it may be seen that the quadratic forms as in the statement of 
Corollary \ref{cor:qf} and their scalar multiples together form a set of full measure in the space
of all binary quadratic forms;  our conclusion shows that for each of these forms the set of
values over the set of pairs $\{(a,b)\mid a,b \in \frak G\}$ is a dense subset of $\C$. 
}
\end{remark}

\noindent {\it Acknowledgements} The first-named author would like to thank the 
Institut de Math\'ematiques de Luminy, Marseille, France and Centre National 
de Recherche Scientifique, France for hospitality and support while this work was
done. 


\vskip15mm
\begin{tabular}{ll} 
S.G. Dani & Arnaldo Nogueira\\ 
School of Mathematics &Institut de Math\'ematiques de Luminy\\
Tata Institute of Fundamental Research& 163, avenue de Luminy, Case 907\\
Homi Bhabha Road &  13288  Marseille  Cedex 9, France\\
Mumbai 400 005, India & \\[4mm]
dani@math.tifr.res.in& nogueira@iml.univ-mrs.fr\\
\end{tabular}

\end{document}